\documentclass[a4paper, 10pt]{article}

\usepackage[T1]{fontenc} 
\usepackage[utf8]{inputenc}
\usepackage{lmodern}
\usepackage[UKenglish]{babel} 
\usepackage{amsmath, amssymb, amsthm, amsfonts}
\usepackage{changepage}
\usepackage[inline, shortlabels]{enumitem}
\usepackage{framed}

\usepackage[a4paper, left=4cm, right=4cm]{geometry}

\usepackage[textsize=scriptsize, bordercolor=lightgray, linecolor=lightgray, backgroundcolor=white]{todonotes}

\usepackage{csquotes}
\usepackage[backend=biber, style=alphabetic, doi=false, isbn=false, url=false]{biblatex}
\renewbibmacro{in:}{}
\usepackage{bookmark}
\hypersetup{colorlinks=true, linkcolor=blue, citecolor=blue}

\usepackage{stmaryrd} 
\usepackage{mathabx} 

\theoremstyle{plain}
\newtheorem*{setting*}{Setting}
\newtheorem*{theorem*}{Theorem}
\newtheorem*{conjecture*}{Conjecture}
\newtheorem*{fact*}{Fact}
\newtheorem{fact}{Fact}
\newtheorem{question}{Question}
\newtheorem{theorem}{Theorem}
\newtheorem{corollary}{Corollary}
\newtheorem{lemma}{Lemma}

\newtheorem*{definition*}{Definition}
\newtheorem*{notation}{Notation}
\newtheorem{proposition}{Proposition}[theorem]

\theoremstyle{definition}
\newtheorem*{remark}{Remark}
\newtheorem*{remarks}{Remarks}
\newtheorem{step}{Step}[proposition]
\newenvironment{proofclaim}{\begin{proof}}{\end{proof}}

\newcommand{\cR}{\mathcal{R}}

\newcommand{\fa}{\mathfrak{a}}
\newcommand{\fb}{\mathfrak{b}}
\newcommand{\fc}{\mathfrak{c}}
\newcommand{\fg}{\mathfrak{g}}
\newcommand{\fG}{\mathfrak{G}}
\newcommand{\fh}{\mathfrak{h}}
\newcommand{\ft}{\mathfrak{t}}
\newcommand{\fu}{\mathfrak{u}}

\newcommand{\fsl}{\mathfrak{sl}}
\newcommand{\fso}{\mathfrak{so}}
\newcommand{\fga}{\mathfrak{ga}}
\newcommand{\fgl}{\mathfrak{gl}}

\newcommand{\GL}{\mathrm{GL}}
\newcommand{\PGL}{\mathrm{PGL}}

\newcommand{\bC}{\mathbb{C}}
\newcommand{\bF}{\mathbb{F}}
\newcommand{\bK}{\mathbb{K}}
\newcommand{\bL}{\mathbb{L}}
\newcommand{\bN}{\mathbb{N}}
\newcommand{\bQ}{\mathbb{Q}}
\newcommand{\bZ}{\mathbb{Z}}

\newcommand{\ad}{\operatorname{ad}}
\newcommand{\DefEnd}{\operatorname{DefEnd}}
\newcommand{\Id}{\operatorname{Id}}
\newcommand{\im}{\operatorname{im}}
\newcommand{\Sym}{\operatorname{Sym}}

\renewcommand{\d}{\partial}

\newcommand{\CZ}{\mathrm{CZ}}

\usepackage{comment}
\excludecomment{MT}

\title{Simple Lie rings of Morley rank $4$\\\Large (The Spanish Inquisition)}
\author{Adrien Deloro and Jules Tindzogho Ntsiri}

\begin{document}

\maketitle

{\footnotesize
\begin{adjustwidth}{6.8cm}{1.6cm}
\begin{itemize}[label=--, itemsep=-2pt]
\itshape
\item
Ici, Chimère ! Arrête-toi !
\item
Non ! Jamais ! 
\end{itemize}
\end{adjustwidth}

\begin{center}
\begin{adjustwidth}{2cm}{2cm}
\noindent\textbf{Abstract.}
We prove the Lie ring equivalent of the Cherlin-Zilber conjecture \textbullet~in characteristic $0$, for any rank and \textbullet~in characteristic $\neq 2, 3$, for rank $\leq 4$. Both are open in the group case.
\end{adjustwidth}
\end{center}
}

\begin{center}
\S~\ref{S:introduction}.~Introduction%
\quad---\quad%
\S~\ref{S:general}.~Generalities%
\quad---\quad%
\S~\ref{S:proofs}.~Proofs%
\quad---\quad%
\S~\ref{S:questions}.~Questions
\end{center}

\section{Introduction}\label{S:introduction}

For first glance, our main result is here; an explanation will follow.
Just consider Morley rank 
as a form of abstract dimension not related to any linear structure.

\begin{theorem*}
Let $\fg$ be a simple Lie ring of finite Morley rank.
\begin{itemize}
\item
If the characteristic is $0$, then $\fg$ is a Lie algebra over a definable field (and one of Lie-Chevalley type).
\item
If the characteristic is $\neq 2, 3$, then $\fg$ does not have Morley rank $4$.
\end{itemize}
\end{theorem*}

The first ad should be considered folklore; as a matter of fact, it was announced without a proof in \cite[Theorem~2]{ZUncountably}. Later \cite{NNonassociative} did not explicitly prove nor state it, but he certainly knew it. This is classical model theoretic-algebra, and 
not the core of the present work. Like the unpublished \cite{Raleph} we shall focus on positive characteristic, where the topic is more algebraic and can be appreciated without knowing mathematical logic. Experts in modular Lie algebras should understand our results, our methods, and most of our final questions.

The second ad is a partial analogue in model theory of \cite[Theorem~2.2]{SLie}. We emphasize that one starts with a Lie ring devoid of any linear structure or gradation, which highly complicates matters: paradoxically, Lie rings bear extremely little geometric information. Therefore our proof does not proceed through Cartan subalgebras, but through Borel subalgebras, with obvious influences from group-theoretic `local analysis'.
Finally we mention that the equivalent for simple \emph{groups} of finite Morley rank is notoriously open and challenging.

Returning to simple Lie rings of finite Morley rank, one conjectures that classifying them would be to the Block-Premet-Strade-Wilson theorem (see~\S~\ref{s:abstractLie}), what the theory of groups of finite Morley rank is to the classification of the finite simple groups: an abstract, simpler sketch (see~\S~\ref{s:modeltheory}). This will deserve future work.


\subsection{Abstract Lie rings}\label{s:abstractLie}

\begin{MT}
\paragraph{Lie and Chevalley.}

The name `Lie' points at two distinct notions. First, a certain type of compatible higher structure on natural group-theoretic objects such as $(\GL_n(\bC); \cdot)$. Second, a certain class of non-associative algebraic structures such as $(\fgl_n(\bC); [\cdot, \cdot])$.
Relations between the two sets of ideas were first discovered by Lie in the intended model of differential geometry, and then pushed to general grounds by Chevalley. How much Lie theory is better understood in general terms, with no specific reference to the intended model of geometry, is however \emph{taboo} in first expositions of the topic. Some mathematicians seem unaware of the failure of categoricity and of the legitimacy of non-Archimedean analysis, not to mention positive characteristic fields. The mathematical physicist will righteously object that `Lie groups' in the physical sense of the term are non-functorial entities. Through provocation we are merely stating the obvious: Lie theory and Chevalley theory are not fully identical, and the focus on higher-order structure is an obstacle to developing the latter.
\end{MT}

\paragraph{Lie-Chevalley functors, Lie-Cartan functors.}

\begin{MT}
Familiar objects such as $\fsl_2$ may be taken over various fields. More generally,
\end{MT}
Chevalley's celebrated basis theorem \cite[\S~25 and notably \S~25.4]{HLinear} brought into Lie algebras the idea of functoriality, with all its power and might. 
Let $\Phi$ be a simple root system, and fix a Chevalley basis of $\Phi$; this gives rise to the Lie ring $\fG_\Phi(\bZ)$. ($\fG$ is upper-case fraktur g.)
Now for $R$ an associative ring with identity, let $\fG_\Phi(R) = \fG_\Phi(\bZ) \otimes_\bZ R$; this is functorial in $R$. The construction $\fG_\Phi$ could be called a (simple) Lie-Chevalley functor; our terminology may evolve in the future.
It generates phrases such as `Lie rings of Lie-Chevalley type', for images of Lie-Chevalley functors, in Chevalley's notation $A_n, \dots, G_2$.
(Experts call `Lie-Chevalley type' the \emph{classical} algebras. This clashes with group-theoretic terminology, for which \emph{classical} only means $A_n, B_n, C_n, D_n$, while $E_6, E_7, E_8, F_4, G_2$ are refered to as \emph{exceptional}.)

The Lie-Chevalley functors, familiar in algebraic geometry, contrast sharply with other constructions yielding simple Lie algebras, usually called `of Cartan type': Witt $W(m; \underline{1})$, special $S(m; \underline{1})^{(1)}$, Hamiltonian $H(m; \underline{1})^{(2)}$, contact $K(m; \underline{1})^{(1)}$. We prefer `of Lie-Cartan type' to prevent confusion with Cartan's maximal toral subalgebras.

While infinite-dimensional in general, simple Lie algebras of Lie-Cartan type become finite-dimensional in positive characteristic. Thus finite-dimensional, simple Lie algebras over algebraically closed fields of positive characteristics do not all fall into Chevalley's families.

\paragraph{Finite-dimensional, simple Lie algebras over algebraically closed fields.}

A key fact is the Block-(Premet-Strade)-Wilson \cite{BWClassification, PSClassification} positive answer to the Kostrikin-Shavarevitch conjecture \cite{KSGraded}, with the following consequence.
\begin{fact*}
Let $\fg$ be a finite-dimensional, simple Lie algebra over an algebraically closed field of characteristic $\geq 5$.
Then :
\begin{itemize}
\item
either $\fg$ is of Lie-Chevalley type,
\item
or $\fg$ is of Lie-Cartan type,
\item
or the characteristic is $5$ and $\fg$ is `of Melikian type' \cite{MSimple}.
\end{itemize}
\end{fact*}

One recommends \cite{PSClassification} or \cite[Introduction]{SSimple1}.
No understanding of the proof nor of the result itself is needed here, as the only simple Lie algebra encountered in the present work is $\fsl_2$.

%

\paragraph{Abstract Lie rings.}

A
\emph{Lie ring} is an abelian group equipped with a bi-additive map satisfying Jacobi's identity and called the \emph{bracket}. Subrings and ideals are defined in these terms; so are nilpotence, solubility, and simplicity. A \emph{Lie morphism} is a morphism of Lie rings, viz.~an additive morphism preserving the bracket. There is a legitimate tradition of calling Lie rings `Lie algebras over $\bZ$'.

The \emph{characteristic} of a Lie ring is the exponent of the underlying additive group; `characteristic $0$' means torsion-free divisible. Every simple Lie ring has characteristic $0$ or a prime.
Every such is therefore a simple Lie algebra over either $\bQ$ or $\bF_p$; but the dimension is in general infinite.
We reserve the phrase \emph{Lie algebra} for the presence of an interesting field acting compatibly. (For us, interesting means `infinite and definable', in the sense below.)

Simple Lie rings are unclassifiable, unless mathematical logic itself brings some form of `tameness'.

\subsection{Lie rings in model theory}\label{s:modeltheory}

The present work aims at showing that some simple Lie rings provided by model theory are \emph{finite-dimensional} Lie algebras.
We restrict our study to finite Morley rank.

\paragraph{Morley rank (for non-logicians).}

When doing model-theoretic algebra, Morley rank is best construed as a dimension on definable sets; see for instance \cite[chapters 2, 3, 4]{BNGroups}. Alternative discussions are \cite[Introduction]{PStable} or \cite[\S~A.I.2]{ABCSimple}.

Let $(\fg; +, [\cdot])$ be a Lie ring. A subset $X\subseteq \fg^n$ is \emph{definable} if there is a first-order formula, viz.~one using the algebraic operations, elements of $\fg$ as parameters, equations, negations, conjunctions, quantification on elements, and giving $X$ as a solution set. (This is the natural generalisation of \emph{constructible} when losing the Chevalley-Tarski theorem that the constructible class is closed under projections.)
For good measure and up to abusing terminology, allow quotients of definable sets by definable equivalence relations.

Now $\fg$ has \emph{finite Morley rank} if there is a dimension function assigning to each non-empty definable set an integer, and satisfying the following:
\begin{itemize}
\item
for definable $X$, $\dim X \geq n+1$ iff there are infinitely many, pairwise disjoint $Y_i \subseteq X$ with $\dim Y_i \geq n$;
\item
for definable $f\colon X \to Y$ and integer $k$, the set $Y_k = \{y \in Y: \dim f^{-1}(\{y\}) = k\}$ is definable;
\item
as above, if $Y = Y_k$, then $\dim X = \dim Y + k$;
\item
as above, there is an integer $m$ such that $Y_0 = \{y \in Y: \operatorname{card} f^{-1}(\{y\}) \leq m\}$.
\end{itemize}

The definition (actually the Borovik-Poizat axiomatisation of groups of finite Morley rank) was given for self-containedness. In practice, there is a dimension provided by mathematical logic, which obeys algebraically predictable rules. For example, `finite non-empty' is equivalent to `$0$-dimensional'.
If $\fg$ is a $\bK$-Lie algebra and $\bK$ is definable and infinite, then $\dim \fg = \dim \bK \cdot \operatorname{ldim}_\bK (\fg)$, where $\operatorname{ldim}_\bK$ denotes the $\bK$-linear dimension. If $\bK$ is finite but $\fg$ is not, the formula makes no sense as $0 \cdot \infty$ is not defined. If $\bK$ is not definable, the formula makes no sense at all since only definable sets bear a dimension.

\paragraph{Lie rings of finite Morley rank.}
So far they received unsufficient attention. Nesin's seminal \cite{NNonassociative}, from which the first author learnt a lot, and Rosengarten's thesis \cite{Raleph}, are not widely cited.
However Baudisch's groups (non-algebraic, nilpotent groups of Morley rank $2$) are built from Lie rings \cite{BNew}, so interactions between Lie rings and model theory are not to be despised. As a matter of fact, an ambitious model-theoretic study of nilpotent Lie algebras and the \emph{Lazard(-Malcev) correspondence} with nilpotent groups is in preparation \cite{DMRSModel} under neo-stability assumptions. But we focus on model-theoretic algebra of finite Morley rank and the following.

\begin{conjecture*}[{`$\log \CZ$'; implicit in \cite{Raleph}}]
Let $\fg$ be an infinite simple Lie ring of finite Morley rank.
Suppose the characteristic is sufficiently large.
Then $\fg$ is a simple Lie algebra over an algebraically closed field.
\end{conjecture*}

\begin{remarks}\leavevmode
\begin{itemize}
\item
Applying the Block-Premet-Strade-Wilson theorem, one would even know the isomorphism type of $\fg$. 
\item
It remains to explain `sufficiently large'. 
The lazy version is: characteristic $0$, which we establish in Theorem~\ref{t:0}.
The modest version is: $> f(d)$ where $d$ is the Morley rank, and $f$ would be determined while working on the conjecture.
The optimistic version is: $> 3$.
There is an intermediate version in characteristic $> d$, where one would also prove that $\fg$ is of Lie-Chevalley type.
\end{itemize}
\end{remarks}


\begin{MT}
\paragraph{Finite-dimensionality, and Cherlin-Macintyre-Shelah properties.}
The full strength of finite Morley rank may be unnecessary.
Recall the following.

\begin{definition*}[{\cite{WDimensional}}]
A theory $T$ is {\upshape [}fine, integer-valued{\upshape ]} \emph{finite-dimensional} if there is a dimension function $\dim$ from the collection of all interpretable sets in models of $T$ to $\bN\cup\{-\infty\}$, satisfying for a formula $\varphi(x, y)$ and interpretable sets $X$ and $Y$:
\begin{itemize}
\item
Invariance: If $a \equiv a'$ then $\dim(\varphi(x, a)) = \dim(\varphi(x, a'))$.
\item
Algebraicity: $X$ is finite nonempty iff $\dim(X) = 0$, and $\dim(\emptyset) = - \infty$.
\item
Union: $\dim(X \cup Y) = \max\{\dim(X), \dim(Y)\}$.
\item
Fibration: If $f\colon X \to Y$ is an interpretable map such that $\dim(f^{-1}(y)) \geq d$ for all $y \in Y$, then $\dim(X) \geq \dim(Y) + d$.
\end{itemize}
\end{definition*}

This generalises both finite Morley rank and $o$-minimal structures. We turn to fields.

An algebraic structure or theory may have the following property: \emph{every infinite definable skew-field is an algebraically closed field }. We call it \emph{the Cherlin-Macintyre-Shelah property} (\textsc{cms}); history tends to forget Zilber here.

The \textsc{cms} holds in universes of finite Morley rank, but does not in arbitrary finite-dimensional universes, notably because of $o$-irreducibility. One may introduce a weak form, w\textsc{cms}: \emph{every infinite definable skew-field is real closed, algebraically closed, or quaternionic over real closed}.

It is not known whether the w\textsc{cms} holds in finite-dimensional universes.
In any case, it is unclear to us if the present work could be carried to `finite-dimensional with w\textsc{cms}' (and one should expect more simple $3$-dimensional rings, at the very least $\fso_3(\cR)$ for real closed $\cR$). But `finite-dimensional with \textsc{cms}' should certainly suffice. Indeed, our focus on connected subgroups could entirely avoid the use of the `descending chain condition' (\textsc{dcc}) on definable subgroups.
\end{MT}

\paragraph{Relations to groups and to Lie algebras.}

\begin{enumerate}
\item
Relation to the Cherlin-Zilber conjecture. The conjecture is a clear analogue of the Cherlin-Zilber conjecture on infinite simple groups of finite Morley rank \cite{ABCSimple}.
We do \emph{not} believe there is a general method for obtaining a Lie ring and an adjoint action from an abstract group of finite Morley rank. Conversely, we do \emph{not} believe there is a general method for retrieving a definable group acting on an abstract Lie ring of finite Morley rank. In short we do \emph{not} believe in a `Lie-Chevalley correspondence' at this level of generality, so the two lines of thought should remain independent.
We do not know whether adding a Zariski geometry \cite{HZZariski} would favorably change the landscape but this is worth asking. See final questions in \S~\ref{S:questions}.
\item
Relation to the Block-Premet-Strade-Wilson theorem.
Work on the Cherlin-Zilber conjecture is supposed to provide a simplified skeleton of the classification of the finite simple groups. Though the latter is settled but not the former, short versions are always interesting to have.
Likewise, we ask experts whether work on `$\log\CZ$' could give a blueprint of the Block-Premet-Strade-Wilson classification \cite{SSimple1, SSimple2, SSimple3}.
\end{enumerate}


\paragraph{Rosengarten's work.}

The `$\log \CZ$' conjecture implicitly motivated Rosengarten's early and unpublished (and in our opinion, underrated) study \cite{Raleph} of Lie rings of Morley rank $\leq 3$. Directed by Cherlin, Rosengarten proved the following in characteristic $\neq 2, 3$:
\begin{itemize}
\item
a connected Lie ring of Morley rank $1$ is abelian \cite[Corollary~4.1.1]{Raleph};
\item
a connected Lie ring of Morley rank $2$ is soluble; if non-nilpotent, it covers $\fga_1(\bK)$ \cite[Theorem~4.2.1]{Raleph};
\item
a connected, non-soluble Lie ring of Morley rank $3$ covers $\fsl_2(\bK)$ \cite[Theorem~4.4.1]{Raleph}.
\end{itemize}

These statements are discussed further in \S\S~\ref{s:1}--\ref{s:3}.

\paragraph{Our results.}
We prove the following.

\begin{theorem*}
Let $\fg$ be a simple Lie ring of finite Morley rank.
\begin{itemize}
\item
If the characteristic is $0$, then $\fg$ is a Lie algebra over a definable field (and one of Lie-Chevalley type).
\item
If the characteristic is $\neq 2, 3$, then $\fg$ does not have Morley rank $4$.
\end{itemize}
\end{theorem*}

Group theorists will deem these results merely inspirational insofar as they should have no impact on the Cherlin-Zilber conjecture, other than a renewed challenge.
However seen from the legitimate Lie-theoretic perspective 
we do believe the connection with the Block-Premet-Strade-Wilson theorem deserves serious investigations in model-theoretic algebra. 

\section{Notation, facts, lemmas}\label{S:general}

This section contains general notation and terminology (\S~\ref{s:notation}), some facts from the theory of groups of finite Morley rank (\S~\ref{s:facts}), and then a couple of lemmas on Lie rings of finite Morley rank (\S\ref{s:lemmas}).

\subsection{Notation; Borel and Cartan subrings}\label{s:notation}

\begin{itemize}
\item
We occasionally combine bracket and product notations, as in Jacobi's identity:
\[[a, bc] = [ab, c] + [b, ac].\]
\item
To distinguish between 
mere additive subgroups, subrings, and ideals, we reserve:
\[\leq\text{ for subgroups}, \qquad \sqsubseteq\text{ for subrings}, \qquad \trianglelefteq\text{ for ideals}.\]
Throughout, $\oplus$ is used at the group-theoretic level.
\item
If $\fg$ is a Lie ring and $x \in \fg$, we denote by $\ad_x$ the derivation $y \mapsto [x, y]$.
\item
We use capital $B$ for `bracket', not for `Borel'. If $\fg$ is a Lie ring and $x \in \fg$, we let $B_x = [x, \fg] = \im \ad_x \leq \fg$, a mere \emph{subgroup}. If $\fg$ is connected (see \S~\ref{s:facts}), so is $B_x$. (Notation generalised below.)
\item
For Borel subrings we favour $\fb \sqsubseteq \fg$. By definition, a Borel subring is one definable, connected, soluble, and maximal as such.
\item
We use capital $C$ for `centraliser', not for `Cartan'. If $\fg$ is a Lie ring with connected components (see \S~\ref{s:facts}) and $x \in \fg$, we let $C_x = C_\fg^\circ (x) = \{y \in \fg: [x, y] = 0\}^\circ\sqsubseteq \fg$, a \emph{subring}.
(Notation generalised below.)
\item
For Cartan subrings we favour $\fc \sqsubseteq \fg$. By definition, a Cartan subring is one definable, connected, nilpotent, and quasi-self-normalising, viz.~of finite index in its normaliser.
\item
For $n$ an integer and $x \in \fg$ we let:
\[B_x^n = \im \ad_x^n \leq \fg \quad\text{and}\quad C_x^n = \ker^\circ \ad_x^n \leq \fg.\]
Be careful that for $n > 1$, $C_x^n$ need not be a subring.
If $\fg$ is connected, these are series of definable, connected sub\emph{groups} (however see Lemma~\ref{l:Cn}); the latter is increasing while the former is decreasing. When one series stabilises, so does the other.
\item
If $x \in \fg$ and $k$ is an integer (modulo the characteristic of $\fg$), we let $E_k(x) = \{y \in \fg: [x, y] = ky\}^\circ = \ker^\circ (\ad_x - k \Id)$.
It is a subgroup, and $[E_k(x), E_\ell(x)] \leq E_{k+\ell}(x)$. Often $x$ will remain implicit. For fixed $x$, the various $E_k$ are in direct sum, being well-understood that indices are integers modulo the characteristic.
\end{itemize}

\subsection{Facts from the theory of groups of finite Morley rank}\label{s:facts}

We need only basic results, used with no reference.
\begin{fact*}
Work inside a group of finite Morley rank.
\begin{itemize}
\item
Every descending chain of definable subgroups is stationary \cite[\S~5.1]{BNGroups}.
\item
Every definable subgroup $A$ has a connected component $A^\circ$, viz.~a smallest definable subgroup of finite index; as a subgroup it is characteristic in $A$ \cite[\S~5.2]{BNGroups}. In case $\fg$ is a Lie ring, then $\fg^\circ$ is an ideal of $\fg$ \cite[Lemma~1]{NNonassociative}. In particular, if $\fg$ is simple and infinite then it is connected.
\item
The image of a definable, connected group under a definable group morphism remains definable and connected.
\item
The sum of two definable, connected subgroups remains definable and connected; the sum of infinitely many such reduces to a finite sum. 
\item
(Strong) Cherlin-Macintyre-Shelah property: infinite definable skew-fields are algebraically closed fields \cite[Theorem~8.10]{BNGroups}.
\end{itemize}
\end{fact*}

\begin{remark}[more on fields]
Augmented fields abound in model theory, even in universes of finite Morley rank. The existence of fields of finite Morley rank with non-minimal multiplicative group (at least in characteristic $0$, \cite{BHMWBoese}) complicates considerably the study of groups of finite Morley rank. They play no role here. However, fields of finite Morley rank with non-minimal additive groups (which exist only in positive characteristic, \cite{BMZRed}) make matters non-trivial. See Steps~\ref{p:3:nonilpotent:st:Bb}, \ref{p:4:nonilpotent:st:contradiction}.
\end{remark}

We turn to Lie rings of finite Morley rank. The following are straightforward adaptations of their group counterparts. More details can be found in \cite[\S~3.3]{Raleph}.

\begin{fact*}
Let $\fg$ be a Lie ring of finite Morley rank.
\begin{itemize}
\item
For $x \in \fg$, one has $\dim \fg = \dim B_x + \dim C_x$.

In particular, if $\fg$ is connected, soluble and infinite, then all centralisers are infinite; otherwise $\fg = B_x \leq \fg'$, a contradiction.
Without solubility there are some challenging open problems; see the questions in \S~\ref{S:questions}.
\item
If $\fg$ is connected, then the ideal $\fg' \trianglelefteq \fg$ generated by commutators is definable and connected.

(This is actually trivial, as abelianity completely bypasses the Chevalley-Zilber `indecomposable generation' lemma of \cite[\S~5.4]{BNGroups}. We take it as an indication that all the strength of finite Morley rank may not be required. See the final questions in \S~\ref{S:questions}.)
\item
If $\fg$ is connected and the centre $Z(\fg)$ is finite, then $\fg/Z(\fg)$ is centreless \cite[Lemma~6.1]{BNGroups}.
\item
If $\fg$ is non-abelian and definably simple (viz.~has no definable ideals other than $\{0\}$ and $\fg$), then $\fg$ is simple.

(To our surprise, this is absent from Rosengarten's work and only implicit in Nesin's. However the argument is present in \cite[top of p.~130]{NNonassociative}.
)
\item
If $\fg$ is nilpotent and infinite, then $Z^\circ(\fg)$ is infinite; every definable, infinite subring meets $Z^\circ(\fg)$, and the normaliser condition holds \cite[Lemmas~6.2 and 6.3]{BNGroups}.
\item
Borel subrings are quasi-self-normalising. (Relies on Fact~\ref{f:1} below.)
\end{itemize}
\end{fact*}

We move to modules, linearisation principles, and their consequences. For $V$ an abelian group, we let $\DefEnd(V)$ be the Lie ring of \emph{definable} endomorphisms. Be careful that $\DefEnd(V)$ itself need not be definable. The bracket is the usual one: $\llbracket f_1, f_2 \rrbracket = f_1\circ f_2 - f_2\circ f_1$.

\begin{definition*}
A \emph{definable Lie ring action} is a triple $(\fg, V, \cdot)$ where:
\begin{itemize}
\item
$\fg$ is a definable Lie ring,
\item
$V$ is a definable, connected abelian group,
\item
$\cdot\colon \fg \times V \to V$ is a definable map inducing a Lie morphism $\rho\colon \fg \to \DefEnd(V)$.
\end{itemize}
When the context is clear we simply say \emph{definable module}.

A definable module is \emph{faithful} if $\ker \rho = \{0\}$, and \emph{$\fg$-irreducible} if $V$ has no definable, connected, infinite, proper $0 < W < V$ which is $\fg$-invariant. (Also called `$\fg$-minimal'.)
\end{definition*}

The following fact extends Zilber's original `definabilisation' of the Schur covariance field.

\begin{fact*}[{consequence of \cite{DZilber}}]
Let $(\fg, V)$ be a definable, faithful, $\fg$-irreducible module in a finite-dimensional universe. Suppose that $\fg$ and $C_{\DefEnd(V)}(\fg)$ are infinite. Then there is a definable field $\bK$ such that $V$ is a finite-dimensional $\bK$-vector space and $\fg$ acts $\bK$-linearly.
\end{fact*}

%
%
%

\subsection{General lemmas}\label{s:lemmas}

We first discuss iterated centralisers.

\begin{lemma}\label{l:Cn}
Let $\fg$ be a Lie ring of finite Morley rank and $h \in \fg$. Let $n$ be an integer such that $C_h^n$ is abelian. Then $C_h^{n+1}$ is a subring.
\end{lemma}
\begin{proof}
Let $\delta = \ad_h$, so that $C_n = \ker^\circ \delta^n$. Let $x, y \in C_{n+1}$. Then:
\[\delta^{n+1}([x, y]) = \sum_{k = 0}^{n+1} \binom{n}{k} [\delta^k x, \delta^{n+1-k} y] = \sum_{k = 1}^n \binom{n}{k} [\delta^k x, \delta^{n+1 - k} y].\]
Now for $1 \leq k \leq n$, one has $\delta^k x \in C_{n+1 - k} \leq C_n$ and $\delta^{n+1 - k} y \in C_k \leq C_n$. By abelianity, all brackets vanish, so $[x, y] \in \ker \delta^{n+1}$. Since $[C_h^{n+1}, C_h^{n+1}]$ is connected, it is contained in $\ker^\circ \delta^{n+1} = C_h^{n+1}$.
\end{proof}

There is a form of converse.

\begin{lemma}[{\cite[Lemma~3.2.6]{Raleph}}]\label{l:divisors}
Let $\fg$ be a Lie ring with no additive $2$-torsion and $a \in \fg$ be such that $C_a^2 = \fg$. Then $B_a$ is an abelian subring of $C_a$, and actually an ideal of $C_a$.
\end{lemma}
\begin{proof}
By assumption, $B_a \leq C_a$. Now let $b_1, b_2 \in B_a$, say $b_i = [a, x_i]$. Then:
\begin{align*}
[b_1, b_2] & = [b_1, [a, x_2]] = [b_1 a, x_2] + [a, b_1 x_2]\\
& = [a, b_1 x_2] = [a, [a, x_1] x_2]\\
& = [a, [ax_2, x_1] + [a, x_1 x_2]]\\
& = [a, [b_2, x_1]] + [a, \underbrace{[a, x_1 x_2]]}_{\in B_a \leq C_a}\\
& = [a b_2, x_1] + [b_2, ax_1]\\
& = [b_2, b_1].
\end{align*}
Since there is no $2$-torsion, we find $[b_1, b_2] = 0$, as wished.

Last, if $b = [a, x] \in B_a$ and $c \in C_a$, then $[b, c] = [[a, x], c] = [ac, x] + [a, xc] = [a, xc] \in B_a$, meaning $B_a \trianglelefteq C_a$.
\end{proof}

The following plays a key role when removing nilpotent Borel subrings in Propositions~\ref{p:3:nonilpotent} and~\ref{p:4:nonilpotent}.

\begin{lemma}[{\cite[Claim~4.3.2]{Raleph}}]\label{l:Rosengarten}
Let $\fg$ be a Lie ring with no additive $2$-torsion and $\fa \sqsubseteq \fg$ be an abelian subring. For $n \geq 1$ and $x \in \fg$, let $A_n = [\fa, \dots, [\fa, x]{\scriptscriptstyle \dots}] \leq \fg$, where $\fa$ appears $n$ times.
Then $[A_n, A_n] \leq [\fa, \fg]$.
\end{lemma}
\begin{proof}
Compute modulo $[\fa, \fg]$.
For any $a \in \fa$ and $x, y \in \fg$ one has:
\[0 \equiv [a, [x, y]] = [ax, y] + [x, ay],\]
whence $[ax, y] \equiv - [x, ay]$. Let $f$ be a product of $n$ $\fa$-derivations, say $f = \ad_{a_1}\circ \cdots\circ \ad_{a_n}$ with the $a_i$ all in $\fa$. By abelianity, the order does not matter. Then $[f(x), y] \equiv (-1)^n [x, f(y)]$. If $f'$ is another product, of say $n'$ $\fa$-derivations, then keeping abelianity in mind:
\begin{align*}
[f(x), f'(x)] & \equiv (-1)^n [x, ff'(x)] \equiv (-1)^n [x, f'f(x)]\\
& \equiv (-1)^{n+n'} [f'(x), f(x)] \equiv (-1)^{n+n'+1} [f(x), f'(x)].
\end{align*}
Hence if $n + n'$ is even then $2 [f(x), f'(x)] \equiv 0$ and $[f(x), f'(x)] \equiv 0$. This is the case when $n' = n$, implying $[A_n, A_n] \leq [\fa, \fg]$.
\end{proof}

The following extremely useful principle is simply additivity of $\dim$.

\begin{lemma}[`lifting the eigenspace']\label{l:liftingEk}
Let $\fg$ be a Lie ring of finite Morley rank and $X$ be a subquotient module, viz.~$X = V/V'$ where $V' \leq V \leq \fg$ are $\fg$-submodules. If for some $h \in \fg$ and some integer $k$ modulo the characteristic, $E_k^X(h) = \{x \in X: [h, x] = x\}^\circ$ is non-trivial, then $E_k^\fg(h) \leq \fg$ is non-trivial.
\end{lemma}
\begin{proof}
Let $\varphi(x) = [h, x] - kx$, which stabilises both $V$ and $V'$. Let $\pi\colon V \to X$ be the quotient map and $\overline{\varphi}\colon X \to X$ be induced by $\varphi$. Notice $E_k^X = \ker^\circ \overline{\varphi}$.
Let $W = (\pi^{-1}(\ker^\circ \overline{\varphi}))^\circ > V'$. Then $\varphi(W) \leq V'$, so $\ker^\circ \varphi > 0$.
\end{proof}

We now turn to the action of a definable subring on the quotient group.

\begin{lemma}\label{l:preHrushovski}
Let $\fg$ be a simple Lie ring of finite Morley rank and $\fh \sqsubset \fg$ be a definable, connected, proper subring. Then $I = C_\fh^\circ(\fg/\fh)$ is nilpotent.
\end{lemma}
\begin{proof}
Of course $I \trianglelefteq \fh$. Let $I^{[n]}$ be the descending nilpotence series, so that $I^{[n]} \trianglelefteq \fh$.

We claim that for $n \geq 0$, one has $[I^{[n+1]}, \fg] \leq I^{[n]}$. Indeed when $n = 0$, one has:
\[
[I^{[1]}, \fg] = [I I, \fg] \leq [I, I\fg] \leq [I, \fh] \leq I,
\]
and then by induction:
\[[I^{[n+2]}, \fg] = [I I^{[n+1]}, \fg] \leq [I, I^{[n1+]}\fg] + [I^{[n+1]}, I\fg] \leq [I, I^{[n]}] + [I^{[n+1]}, \fh] \leq I^{[n+1]}.\]

By the descending chain condition, there is $n$ such that $I^{[n+1]} = I^{[n]}$. Then $[I^{[n]}, \fg] \leq I^{[n]}$ so $I^{[n]} \trianglelefteq \fg$. Since $I^{[n]} < \fg$, we get $I^{[n]} = 0$ by simplicity.
\end{proof}

To our own surprise we could easily obtain a self-normalisation criterion.

\begin{lemma}\label{l:selfnormalisation}
Let $\fg$ be a definable, connected Lie ring and $\fh \sqsubset \fg$ be a definable, connected subring of codimension $1$. Then either $\fh \triangleleft \fg$ is an ideal, or $N_\fg(\fh) = \fh$.
\end{lemma}
\begin{proof}
Let $I = C_\fg(\fg/\fh)$. If $I = \fh$ then $\fh$ is an ideal of $\fg$; so suppose not.  Then $\fh/I$ acts non-trivially on $\fg/\fh$. Linearising in dimension $1$ (see \S~\ref{s:1}), the action is scalar; hence free.

Let $h \in \fh \setminus I$. Let $n \in N_\fg(\fh)$. Then $[h, n] = - [n, h] \in \fh$, so the non-trivial, scalar action of $h$ on $\fg/\fh$ sends the image of $n$ to $0$. By freeness, $n \in \fh$ already.
\end{proof}

Last is an important observation on definable derivations.

\begin{lemma}\label{l:noderivations}
An infinite field of finite Morley rank has no non-trivial definable derivation.
\end{lemma}
\begin{proof}
Recall that the Cherlin-Macintyre-Shelah property holds: $\bK$, or any infinite definable subfield, is algebraically closed.

Let $\delta \colon \bK \to \bK$ be a definable derivation, and $\bK_0$ be the field of constants. Suppose the characteristic is not $2$. Let $x \in \bK_0$, and let $y \in \bK$ such that $y^2 = x$. Then $0 = \delta(x) = 2 y \delta(y)$, so $y \in \bK_0$. So $\bK_0$ is closed under taking square roots: hence infinite. In characteristic $2$, it is however closed under taking cubic roots: hence infinite as well.
So in either case, $\bK_0$ is infinite and definable: this proves $\bK_0 = \bK$, viz.~$\delta = 0$.
\end{proof}

The same holds in $o$-minimal universes, but seems unclear in the general finite-dimensional case. It however appears like a desirable property for a more systematic investigation of Lie rings in model theory.

\section{The proofs}\label{S:proofs}

\S~\ref{s:0} handles characteristic $0$; Theorem~\ref{t:0} will not surprise people familiar with Zilber's early work. Turning to positive characteristic, \S\S~\ref{s:1}, \ref{s:2}, \ref{s:3}, \ref{s:4} deal with increasing dimensions, and theorems are numbered accordingly. Since we do not reprove Rosengarten's rank $1$ analysis, we call it \emph{Fact~\ref{f:1}}; Corollary~\ref{c:1} should be folklore. Theorems~\ref{t:2} is already in \cite{Raleph}, but the too partial (though important) Corollary~\ref{c:Lie2} is new. Theorem~\ref{t:3} too is already in \cite{Raleph}; its simple Corollary~\ref{c:3} is not. Theorem~\ref{t:4} is definitely new.

\setcounter{subsection}{-1}
\subsection{Characteristic 0}\label{s:0}

The following entirely settles matters in characteristic $0$. It is essentially Zilber's definability of the Schur covariance field and was known to Zilber \cite[Theorem~2]{ZUncountably}. We learnt about this paper from Belegradek, whom we quote (personal correspondence).
\begin{quote}
\itshape\small
This is one-page paper containing only formulations of results and some hints of proofs. It was submitted in 1979 but published only in 1982.
The rule of that journaI was that they could publish a short paper containing no proofs in case if the author also submitted a manuscript with detailed proofs; the referee wrote a report based on the short paper and manuscript together.
I don't know whether that manuscript is still available.
But I think it makes sense for people working on groups of finite Morley rank to rediscover the proofs of the paper, and to advertise the paper (which seems not well-known even though its author is the famous Zilber; maybe this is because that second rate Russian journal was not translated into English, and
the paper was not reviewed in MR).
\end{quote}

Although we found no further trace of the topic in the literature, the result must also have been known to Nesin \cite{NNonassociative} among others. Why Rosengarten does not discuss it, is a mystery to us.

\setcounter{theorem}{-1}
\begin{theorem}\label{t:0}
Work in a theory of finite Morley rank.
\begin{enumerate}[label=(\roman*)]
\item\label{t:0:i:linearisation}
Let $(\fg, V, \cdot)$ be a definable Lie ring action with $\fg$ infinite.
Suppose it is irreducible and faithful. \emph{Suppose further that $V$ does not have bounded exponent.} Then the configuration is linear: there is a definable field $\bK$ such that $V$ is a $\bK$-vector space and $\fg$ is a $\bK$-Lie subalgebra of $\fgl_\bK(V)$.
\item
If in addition $\fg$ is a Lie algebra over some algebraically closed field $\bL$, then $\bL \simeq \bK$ definably.
\item
If $\fg$ is a simple Lie ring of finite Morley rank \emph{and characteristic $0$}, then there is a definable field $\bK$ such that $\fg$ is the ring of $\bK$-points of one of the Lie-Chevalley functors $A_n, \dots, G_2$.
\end{enumerate}
\end{theorem}

\begin{proof}\leavevmode
\begin{enumerate}[label=(\roman*)]
\item
The action of $\fg$ commutes to $\bZ\Id_V$, so we may linearise (see \S~\ref{s:facts}).
There is a definable field $\bK$ such that $V$ a finite-dimensional vector space over $\bK$ and $\fg \hookrightarrow \fgl_\bK(V)$, a mere embedding of Lie rings. Of course $\operatorname{char} \bK = 0$. It suffices to show that $\fg$ is a vector subspace of $\fgl_\bK(V)$. Indeed $N = \{\lambda \in \bK: (\forall x \in \fg)(\lambda \cdot x \in \fg)\}$ is an infinite, definable subring of $\bK$, hence equal to $\bK$ by the Cherlin-Macintyre-Shelah property.
\item
Let $\VDash$ denote interpretability.
The above and definability of linear structures over a pure field imply $(\bK; +, \cdot) \VDash (\fg, V; +, [\cdot, \cdot], \cdot)$. Now using the adjoint action of a Cartan subalgebra, one has $(\fg; +, [\cdot, \cdot]) \VDash (\bL; +, \cdot)$. It follows $(\bK; +, \cdot) \VDash (\bL; +, \cdot)$ and one may conclude by Poizat's monosomy theorem \cite[Theorem~4.15]{PStable} that $\bK \simeq \bL$, even $\bK$-definably.
\item
Let $\fg$ act on $V = (\fg; +)$ by the adjoint representation: the module is faithful and irreducible. By~\ref{t:0:i:linearisation}, $\fg$ is a Lie algebra over some algebraically closed field $\bK$. It remains simple as such. Since $\bK$ has characteristic $0$, the isomorphism type of $\fg$ is given by the classification into Chevalley families.
\qedhere
\end{enumerate}
\end{proof}

Nothing similar is possible in positive characteristic---and one should also expect objects of Lie-Cartan type.

\subsection{Dimension 1}\label{s:1}

We state Rosengarten's `minimal' result (which we do not reprove), and then a number of consequences in Corollary~\ref{c:1}.

\begin{fact}[{Rosengarten; \cite[Theorem~4.1.1]{Raleph}}]\label{f:1}
Let $\fg$ be an infinite Lie ring of finite Morley rank and characteristic $\neq 2, 3$.
Then $\fg$ contains an infinite, definable, connected, abelian subring.
\end{fact}

This is the analogue of Reineke's theorem on `definably minimal' groups \cite{RMinimale}.
The existing proof is non-trivial and relies on results in the finite case. (Moreover the problem is open in characteristic $3$).
See the section on open questions, \S~\ref{S:questions}.

We prove a couple of consequences, some already mentioned in \S~\ref{S:general}. In general they will be used without mention.

\begin{corollary}\label{c:1}
Let $\fg \neq 0$ be a connected Lie ring of finite Morley rank and characteristic $\neq 2, 3$.
\begin{enumerate}[label=(\roman*)]
\item
Borel subrings are infinite and quasi-self-normalising.
\item\label{c:1:i:defmin}
Let $V$ be a definable, faithful  module. Suppose that $V$ is definably minimal as a group (viz.~$\{0\}$-irreducible). Then $\dim \fg \leq \dim V$. 
\item
Let $V$ be a faithful $1$-dimensional $\fg$-module. Then there is a field $\bK$ with $V \simeq \bK_+$ and $\fg \simeq \bK \Id_V$.
In particular, there is $h \in \fg$ acting like $\Id_V$.
\end{enumerate}
\end{corollary}
We call the last ad `linearising in dimension $1$'.
\begin{proof}\leavevmode
\begin{enumerate}[label=(\roman*)]
\item
By Fact~\ref{f:1}, there exist infinite definable abelian subrings, so Borel subrings are non-trivial. Let $\fb \sqsubseteq \fg$ be one such. If $N_\fg(\fb)/\fb$ is infinite, then it contains an infinite definable abelian subring $\overline{\fa}$. Lifting and taking the connected component, the soluble subring $\pi^{-1}(\overline{\fa})^\circ > \fb$ contradicts the definition of $\fb$.
\item
Let $V$ be a faithful $\fg$-module which is $\{0\}$-irreducible.
Let $v \in V\setminus\{0\}$ and $\fh = C_\fg^\circ(v)$. Suppose that $\fh$ is infinite. By Fact~\ref{f:1}, it contains an infinite, definable, connected, abelian $\fa \sqsubseteq \fh$. Now $V$ remains faithful as an $\fa$-module, and is $\fa$-irreducible. Linearising the abelian action (see end of \S~\ref{s:facts}), $V$ is a vector space over a definable field $\bK$, and the action of $\fa$ is by scalars, hence free. But $\fa$ acts trivially on $v \neq 0$: a contradiction. So $\fh$ is finite. In particular, $\dim (\fg \cdot v) = \dim \fg \leq \dim V$.
\item
Let $V$ be a $1$-dimensional, faithful $\fg$-module. By~\ref{c:1:i:defmin}, $\dim \fg = 1$. By Fact~\ref{f:1}, $\fg$ is abelian and we may linearise. Clearly $V \simeq \bK_+$ and $\fg \simeq \bK \Id_V$.
\qedhere
\end{enumerate}
\end{proof}


\subsection{Dimension 2}\label{s:2}

We reprove Rosengarten's dimension $2$ analysis (Theorem~\ref{t:2}), and derive an important, though too partial, form of Lie's theorem \emph{in dimension $2$} (Corollary~\ref{c:Lie2}).
The following is already in \cite{Raleph}; our proof is as short but more conceptual.

\setcounter{theorem}{1}
\begin{theorem}[{Rosengarten; \cite[Theorem~4.2.1]{Raleph}}]\label{t:2}
Let $\fg$ be a connected Lie ring of Morley rank $2$ and characteristic $\neq 2, 3$. Then:
\begin{enumerate}[label=(\roman*)]
\item
$\fg$ is soluble;
\item
if $\fg$ is non-nilpotent, then:
\begin{itemize}
\item
there is a Cartan subring $\fc \sqsubseteq \fg$ with $\fg = \fc \oplus \fg'$; 
\item
there is $h \in \fc$ with $(\ad_h)_{|\fg'} = \Id_{\fg'}$ (in particular, $Z(\fg) \leq \fc$);
\item
$Z(\fg)$ is finite and there is a definable field $\bK$ with $\fg/Z(\fg) \simeq \fga_1(\bK)$.
\end{itemize}
\end{enumerate}
\end{theorem}
\begin{proof}\leavevmode
\begin{enumerate}[label=(\roman*)]
\item
Let $\fh \sqsubseteq \fg$ be an infinite, definable, connected subring of minimal $\dim$. By Fact~\ref{f:1}, $\fh$ is abelian, so we may assume it is proper; $\dim \fh = 1$. Consider the action of $\fh$ on $V = \fg/\fh$. If it is trivial, then $\fh \triangleleft \fg$ is an ideal, so $\fg$ is soluble and we are done. So we suppose that said action is non-trivial. Linearising in dimension $1$, 
%
there is $h \in \fh$ 
acting
on $V$ as $\Id_V$.
%
Lifting the eigenspace with Lemma\ref{l:liftingEk}, $E_1(h) \neq 0$. Now $\fh$ is abelian so $\fh \leq E_0(h)$.
Since $\dim \fg = 2$, we get $\fg = E_0 \oplus E_1$; finally $[E_1, E_1] \leq E_2 = 0$ so $E_1$ is an ideal of $\fg$ and we are done.
%
\item
Suppose $\fg$ is non-nilpotent.
Then $\fg'$ has dimension $1$, so it is abelian by Fact~\ref{f:1}. By non-nilpotence, the action of $\fg$ on $\fg'$ is non-trivial. Linearising in dimension $1$, there is $h \in \fg$ acting on $\fg'$ like $\Id_{\fg'}$; hence $E_1 \neq 0$. Since $\fg$ is soluble, one has $E_0 \neq 0$. So we find $\fg = E_0 \oplus E_1$ and $E_1 = \fg'$. Always by non-nilpotence, $E_1$ is the only $1$-dimensional, connected ideal of $\fg$. It follows that $\fc = E_0$ is abelian, and has finite index in its normaliser: it is a Cartan subring. Clearly $Z(\fg) \leq E_0$. Last, $Z(\fg)$ is finite, so $\fg/Z(\fg)$ is centreless. We now suppose $\fg$ centreless. The action of $E_0$ on $E_1$ is then faithful, so linearisation directly gives the required isomorphism.
%
\qedhere
\end{enumerate}
\end{proof}

\begin{remark}
A finite centre is unavoidable, as the following example shows. Let:
\[\fg = \left\{\begin{pmatrix} x & 0 & 0\\ 0 & \alpha(x) & y\\ 0 & 0 & 0\end{pmatrix}: (x, y) \in \bK^2\right\},\]
for any additive $\alpha\colon \bK_+ \to \bK_+$, such as for instance $x \mapsto x^p$.
Then $Z(\fg)$ is isomorphic to $\ker \alpha$. 
\end{remark}

The following is very unsatisfactory and one should generalise it to $\dim V$ less than the characteristic, with no mention of $\dim \fb$. It plays an important role in the proof of Theorem~\ref{t:4}.

\begin{corollary}\label{c:Lie2}
Work in finite Morley rank.
Let $\fb$ be a non-nilpotent, connected Lie ring of characteristic $\neq 2, 3$. Let $V$ be an irreducible $\fb$-module. Suppose $\dim \fb = \dim V = 2$. Then $\fb'$ centralises $V$.
\end{corollary}
\begin{proof}
By Theorem~\ref{t:2}, $\fa = \fb' > 0$; moreover there is a Cartan subring $\ft < \fb$ such that $\fb = \ft \oplus \fa$ and $\fa = [\ft, \fa]$.

Let $\Delta = \DefEnd(V)$, which is an associative ring with induced Lie bracket $\llbracket f, g\rrbracket = fg - gf$. Let $\rho\colon \fb \to \Delta$, a Lie morphism. We avoid it from notation mostly, with exceptions for enhanced clarity. To prove that $\fa \leq \ker \rho$, we suppose otherwise.
Notice that $C_V^\circ(\fa)$ is $\fb$-invariant, so we may suppose $C_V^\circ(\fa) = 0$.


\renewcommand{\thestep}{\thecorollary.\arabic{step}}
\begin{step}
$V$ is not $\fa$-irreducible.
\end{step}
\begin{proofclaim}
Suppose it is. Linearising the abelian action, $\bK = C_{\DefEnd(V)}(\rho(\fa))$ is a definable field containing $\rho(\fa)$.
We claim that $\ft$ acts by definable derivations of $\bK$. Indeed, let $\lambda \in \bK$ and $a \in \fa$. Also let $h \in \ft$ and $a' = [h, a]$. Then:
\begin{align*}
\llbracket \llbracket \rho(h), \lambda\rrbracket, \rho(a)\rrbracket & = (\rho(h) \lambda - \lambda \rho(h)) \rho(a) - \rho(a) (\rho(h) \lambda - \lambda \rho(h))\\
& = 
\rho(h) \rho(a) \lambda - \lambda \rho(h) \rho(a) - \rho(a) \rho(h) \lambda + \lambda \rho(a) \rho(h)\\
& = \llbracket\rho (h), \rho(a)\rrbracket \lambda - \lambda \llbracket \rho(h), \rho(a)\rrbracket\\
& = \rho([h, a)] \lambda - \lambda \rho([h,a])\\
& = \llbracket \rho(a'), \lambda\rrbracket\\
& = 0,
\end{align*}
proving that $\rho(h)$ does act on $\bK$. The action is clearly that of a derivation. By Lemma~\ref{l:noderivations}, $\rho(\ft)$ must act trivially; in particular $\rho(\ft)$ centralises $\rho(\fa) \leq \bK$.
So $\fa = [\ft, \fa] \leq \ker \rho$, as wanted.
\end{proofclaim}

\begin{step}
$V$ is $\fa$-irreducible.
\end{step}
\begin{proofclaim}
Otherwise let $W < V$ be a $1$-dimensional $\fa$-submodule.
If $\fa$ acts trivially on $W$, then $C_V^\circ(\fa) \neq 0$, a contradiction. So $\fa$ acts non-trivially on $W$.

Linearising, $\bK = C_{\DefEnd(W)}(\fa)$ is a definable field containing $\rho(\fa)$. We fix $h \in \ft$ ad-acting on $\fa$ like $\Id_\fa$, viz.~$[h, a] = a$ for $a \in \fa$.

We claim $hW \neq 0$ and $(W \cap hW)^\circ = 0$, so that $V = W + h W$. Suppose $w_1 = h w_2$ for some $w_1, w_2 \in W$. Then applying an arbitrary $a \in \fa$:
\[a w_1 = a h w_2 = ha w_2 - a w_2,\]
so $h a w_2 = a(w_1 + w_2) \in W$. But if $w_2 \notin C_V(\fa)$, then $\fa w_2 = W$, which proves $h W = W$. So $h$ acts on $W$, and as before on $\bK$ as well. Now it induces a derivation of $\bK \geq \fa$, although it acts non-trivially on $\fa$; against Lemma~\ref{l:noderivations}. So $w_2 \in C_V(\fa)$, which is finite.
So $w_1 = h w_2$ has finitely many solutions, which proves all claims.

We derive a contradiction.
Let $w \in W$. Then there are $w_1, w_2 \in W$ such that $h^2 w = w_1 + hw_2$.
Let $a \in \fa$ act on $W$ like $\Id_W$ (this is the element going to $1_\bK$). Remember $[h, a] = a$; in particular $a h^2 = h^2 a - 2 ha + a$.
Applying to $w$ one finds:
\begin{align*}
a h^2 w & = a w_1 + a h w_2 = w_1 + ha w_2 - a w_2 = w_1 - w_2 + h w_2\\
= (h^2 a - 2 ha + a)w & = h^2 w - 2 h w + w = w_1 + w + h (w_2 - 2w).
\end{align*}
This proves $hw \in W \cap hW$ which is finite, so by connectedness $hW = 0$, a contradiction.
\end{proofclaim}
This completes the proof of Corollary~\ref{c:Lie2}.
\end{proof}
\renewcommand{\thestep}{\theproposition.\arabic{step}}

\subsection{Dimension 3}\label{s:3}

\begin{theorem}[{Rosengarten; \cite[Theorem~4.4.1]{Raleph}}]\label{t:3}
Let $\fg$ be a simple Lie ring of Morley rank $3$ and characteristic $\neq 2, 3$. Then $\fg \simeq \fsl_2(\bK)$.
\end{theorem}

Our three reasons to give a proof are:
\begin{enumerate*}
\item
Rosengarten's own argument was never published;
\item
ours is sufficiently different from Rosengarten's, and more conceptual;
\item
parts of it make sense in the group setting;
\item
some key ideas reappear in the dimension $4$ argument of \S~\ref{s:4}.
\end{enumerate*}

The proof first removes nilpotent Borel subrings (Proposition~\ref{p:3:nonilpotent}), then explicitly identifies $\fg$ (Proposition~\ref{p:3:identification}).

\begin{proposition}\label{p:3:nonilpotent}
No Borel subring is nilpotent.
\end{proposition}
\begin{proof}
Let $\fb \sqsubseteq \fg$ be a nilpotent Borel subring. Let $Z = Z^\circ(\fb) > 0$.
First we prove that $\fb$ acts trivially on subquotients of rank $1$ (Step~\ref{p:3:nonilpotent:st:1modules}), and has rank $1$ itself. Then we find a field $\bL$ with rank $2$, such that $\fb$ embeds properly into $\bL_+$ (Step~\ref{p:3:nonilpotent:st:Bb}). Then we derive a contradiction (Step~\ref{p:3:nonilpotent:st:contradiction}).

\begin{step}\label{p:3:nonilpotent:st:1modules}
Let $X$ be a $1$-dimensional subquotient $\fb$-module. Then $\fb$ centralises $X$.
\end{step}
\begin{proofclaim}
Suppose not. Linearising in dimension $1$, there is $h \in \fb$ acting on $X$ like $\Id_X$. Lifting the eigenspace, $E_1 = E_1(h) \leq \fg$ is non-trivial. It is proper since $E_0 \geq Z > 0$.

For each integer $k \not\equiv 0$ one has $E_k\cap \fb = 0$. Indeed, if $a \in E_k\cap \fb$, then $a = \frac{1}{k} [h, a] \in \fb'$, and for every integer $n$ one has $a \in \fb^{[n]}$, the descending nilpotent series. Hence $a = 0$.

We claim $\dim \fb = 1$. Otherwise $\dim \fb = 2$ and $\fg = \fb \oplus E_1$. But then $E_2 = 0$, so $E_1$ is an abelian subring. Now $Z \leq E_0$ normalises $E_1$. If the action is trivial, then $E_1 \leq C_\fg^\circ(Z) = \fb$, a contradiction. So the action of $Z$ on $E_1$ is non-trivial. Hence we may assume $h \in Z$; in that case, $E_0 = \fb$ normalises $E_1$, which is an ideal of $\fg$. This is a contradiction; hence $\dim \fb = 1$.

It follows that $\fb$ is abelian and $E_0 = \fb$. If $E_2 \neq 0$ then $\fg = E_0 \oplus E_1 \oplus E_2$ and $E_2 \triangleleft \fg$; a contradiction. So $E_2 = 0$ and $E_1$ is an abelian subring. But then $E_0 \oplus E_1$ is a subring properly containing $E_0 = \fb$. By definition of a Borel subring, $\fg = E_0 \oplus E_1$, which contains $E_1$ as an ideal: a contradiction.
\end{proofclaim}

It follows $\dim \fb = 1$. Otherwise $\dim \fb = 2$ and $\dim \fg/\fb = 1$. By Step~\ref{p:3:nonilpotent:st:1modules}, $\fb$ centralises $\fg/\fb$, meaning $\fb \triangleleft \fg$, a contradiction. Hence $\dim \fb = 1$ and in particular $\fb$ is abelian. It follows that for $a \in \fb\setminus\{0\}$ one has $C_a = \fb$.

Let $B_\fb = [\fb, \fg]$.

\begin{step}\label{p:3:nonilpotent:st:Bb}
$\fg = \fb + B_\fb$.
There is a definable field $\bL$ of dimensional $2$ such that $B_\fb \simeq \bL_+$ and $\fb \hookrightarrow \bL \Id_{B_\fb}$.
\end{step}
\begin{proofclaim}
First, $\dim B_\fb = 2$. Indeed, let $a, a' \in \fb$ be non-zero. Clearly $\dim B_a = 2$; by abelianity, $B_a$ is $\fb$-invariant. By Step~\ref{p:3:nonilpotent:st:1modules}, $\fb$ centralises $\fg/B_a$. In particular $B_{a'} \leq B_a$ and equality follows. So $B_\fb = B_a$ for any $a \in \fb\setminus\{0\}$, which proves $\dim B_\fb = 2$.

Suppose $\fb \leq B_\fb$. Then for $a \in \fb \setminus\{0\}$ one has $C_a = \fb \leq B_\fb = B_a$. In particular $B_a^2 < B_a$, and $C_a^2 > C_a$. However, $C_a = \fb$ is abelian, so by Lemma~\ref{l:Cn}, $C_a^2$ is a subring. But $\dim C_a^2 \leq 2 \dim C_a = 2\dim \fb < \dim \fg$, so $C_a^2$ is a definable, connected, proper subring containing $C_a = \fb$. Hence $C_a^2 = \fb = C_a$, a contradiction proving $\fb \not\leq B_\fb$.

$B_\fb$ is $\fb$-irreducible. Otherwise there is a $1$-dimensional $\fb$-submodule $V_1 < B_\fb$. By Step~\ref{p:3:nonilpotent:st:1modules}, $\fb$ centralises $V_1$, so $V_1 \leq C_\fg^\circ(\fb) = \fb$. Always by Step~\ref{p:3:nonilpotent:st:1modules}, $\fb$ centralises $B_\fb/\fb$, so $B_\fb\leq N_\fg^\circ(\fb) = \fb$, a contradiction.

The action of $\fb$ on $B_\fb$ is faithful. Otherwise there is $a \in \fb\setminus\{0\}$ with $B_\fb \leq C_a = \fb$, a contradiction again.

Linearising the abelian action, there is a definable field $\bL$ with $B_\fb \simeq \bL_+$ and $\fb \hookrightarrow \bL \Id_{B_\fb}$.
\end{proofclaim}

\begin{step}\label{p:3:nonilpotent:st:contradiction}
Contradiction.
\end{step}
\begin{proofclaim}
Fix any $a_0 \in \fb \setminus\{0\}$ and $x_0 \in B_\fb\setminus\{0\}$.

Because the action of $\fb$ on $B_\fb$ is by scalars, $\dim [a_0, [\fb, x_0]] = \dim [\fb, x_0] = 1$ (the subgroups ned not be equal).
We claim that $\dim [\fb, [\fb, x_0]] = 1$.
Let $A_2 = [\fb, [\fb, x_0]] \leq B_\fb$. By Lemma~\ref{l:Rosengarten}, $[A_2, A_2] \leq B_\fb$.
If $\dim A_2 = 2$ then dimensions match and $A_2 = B_\fb$. Now $[B_\fb, B_\fb] \leq B_\fb$, so $B_\fb$ is a subring. It also is $\fb$-invariant, but $\fb \not\leq B_\fb$ by Step~\ref{p:3:nonilpotent:st:Bb}. So $B_\fb \triangleleft \fg$, a contradiction. Hence $\dim A_2 = 1$.

Thus $\dim [a_0, [\fb, x_0]] = \dim [\fb, [\fb, x_0]]$; the groups are equal.
So for $a_1, a_2 \in \fb$, there is $a_3 \in \fb$ with $[a_0,[a_3, x_0]] = [a_1, [a_2, x_0]]$. Returning to the field $\bL$, this means $a_0 a_3 = a_1 a_2$ as elements of $\bL$. So the image of $\fb$ in $\bL$ is closed under the map $(a_1, a_2) \mapsto \frac{a_1 a_2}{a_0}$.

Let $R = \{\lambda \in \bL: \lambda \fb \leq \fb\}$, a definable subring. We just proved $\frac{1}{a_0} \fb \leq R$, so $R$ is infinite. By the Cherlin-Macintyre-Shelah property, every infinite definable domain is an algebraically closed field. So $R = \bL$ and $\fb \neq 0$ is an ideal of $\bL$. Hence $\fb = \bL$ but dimensions do not match.
\end{proofclaim}

This eliminates nilpotent Borel subrings in a $3$-dimensional, simple Lie ring: Proposition~\ref{p:3:nonilpotent} is proved.
(The argument will be extended to the $4$-dimensional case in Proposition~\ref{p:4:nonilpotent}.)
\end{proof}

By Theorem~\ref{t:2}, the structure of Borel subrings of $\fg$ is well-understood.

\begin{proposition}\label{p:3:identification}
$\fg \simeq \fsl_2(\bK)$.
\end{proposition}
\begin{proof}
The proof starts by some local analysis, viz.~the study of intersections of Borel subrings (Step~\ref{p:3:identification:st:local}). Then we find the desired weight space decomposition (Step~\ref{p:3:identification:st:weights}), which allows final coordinatisation (Step~\ref{p:3:identification:st:coordinatisation}).

\begin{step}\label{p:3:identification:st:local}
If $\fb_1 \neq \fb_2$ are distinct Borel subrings, then $\fc = (\fb_1\cap \fb_2)^\circ$ is a Cartan subring of both.
\end{step}
\begin{proofclaim}
Clearly $\dim \fc = 1$. By Proposition~\ref{p:3:nonilpotent}, $\fb_1$ and $\fb_2$ are non-nilpotent; it is enough to prove that $\fc$ is normal in neither. Suppose it is, say in $\fb_1$; then $\fc = \fb_1' \leq \fb_2$. Considering $N_\fg^\circ(\fc) < \fg$, the ring $\fc$ is not normal in $\fb_2$, so it is a Cartan subring of $\fb_2$.

Let $\ft$ be a Cartan subring of $\fb_1$. Let $t \in \ft$ act on $\fb'_1 = \fc$ like the identity. Now let $h \in \fc$ act on $\fb_2'$ like the identity.
Notice that $[\ft, h] = \fb'_1 = \fc$, and $[h, \fb'_2] = \fb'_2$, so $B_h \geq \fc + \fb'_2 = \fb_2$. Equality follows since $C_h \geq \fb'_1 > 0$.

Now for $f \in \fb'_2$ one has:
\[[t, f] = [t, [h, f]] = [th, f] + [h, tf] = [h, f] + [h, tf].\]
The right-hand is in $B_h = \fb_2$. Therefore $[t, f] \in \fb_2$, and the right-hand is in $\fb_2'$. Therefore $[t, f] \in \fb_2'$, implying $[h, tf] = [t, f]$. So there remains $[h, f] = 0$, a contradiction.
\end{proofclaim}

\begin{step}\label{p:3:identification:st:weights}
Weight decomposition: there is $h_1 \in \fg$ such that $\fg = E_{-2} \oplus E_0 \oplus E_2$, and $[E_{-2}, E_2] = E_0$. Moreover, for $e \in E_2\setminus\{0\}$ one has $C_\fg(e) = E_2$.
\end{step}
\begin{proofclaim}
Let $\fb \sqsubseteq \fg$ be a Borel suring. By Proposition~\ref{p:3:nonilpotent}, it is non-nilpotent. In particular, by Theorem~\ref{t:2}, there is a definable field such that $\fb/Z(\fb) \simeq \fga_1(\bK)$. Let a temporary $h_1 \in \fb$ act on $\fb'$ like $\Id_{\fb'}$. (It will be rescaled at the end of the proof of the Step.)
Of course $E_0 = E_0(h_1) = C_{h_1} > 0$ and $E_1 \geq \fb' > 0$.
Let $e \in \fb'\setminus\{0\}$.

We claim that $C_e = \fb'$ and $B_e = \fb$. Clearly $C_e \geq \fb'$. If equality does not hold, then $\dim C_e = 2$. Now $C_e$ is a Borel subring distinct from $\fb$ since $e \notin Z(\fb)$. But $C_e \geq \fb'$ which is not Cartan in $\fb$, against Step~\ref{p:3:identification:st:local}. So $C_e = \fb'$ and $\dim B_e = 2$. Now consider the action of $\fb$ on $\fg/\fb$. Linearising in dimension $1$, $\fb'$ acts trivially so $B_e \leq \fb$ and equality holds.

We now claim that $E_{-1} = E_{-1}(h_1) \neq 0$. For suppose $E_{-1} = 0$. The additive morphism $\varphi(x) = [h_1, x] + x$ then has a finite kernel, and is surjective. Therefore:
\[\im \ad_{h_1} \ad_e = \im \ad_e (1 + \ad_{h_1}) = \im \ad_e = B_e = \fb.\]
But $\im \ad_{h_1} \ad_e = [h_1, B_e] = [h_1, \fb] < \fb$, a contradiction proving $E_{-1} \neq 0$.

It follows $\fg = E_{-1} \oplus E_0 \oplus E_1$, and $\fc = E_0$ normalises each term. Also, $\fb = E_0 \oplus E_1$ and $\fb' = E_1$. If $[E_{-1}, E_1] = 0$, then $E_{-1} \oplus E_1$ is a subring, hence a Borel subring. However it intersects $\fb$ in $E_1 = \fb'$, against Step~\ref{p:3:identification:st:local}. So $[E_{-1}, E_1] = E_0$.

There remains to prove that for any $e \in E_1 \setminus\{0\}$, one has $C_\fg(e) = E_1$ (notice the absence of a connected component).
Since $\fg = E_{-1} \oplus E_0 \oplus E_1$, it is enough to prove $C_{E_0}(e) = C_{E_{-1}}(e) = 0$.
\begin{itemize}
\item
Let $h \in E_0$ centralise $e$; then it centralises $E_0$ and $[E_0, e] = E_1$. But also, for $y \in E_{-1}$, one has $[e, y] \in E_0 \leq C_h$, so:
\[[e, [h, y]] = [e h, y] + [h, e y] = 0.\]
Therefore $[h, E_{-1}] \leq C_e = E_1$. However $[h, E_{-1}] \leq E_{-1}$, proving $h \in Z(\fg) = 0$.
\item
Now let $f \in E_{-1}$ centralise $e$. Then it normalises $C_e = E_1$ and $N_\fg^\circ(E_1) = E_0 + E_1$. In particular, $[f, h_1] \in E_0 + E_1$, but $[f, h_1] = - [h_1, f] = f \in E_{-1}$ so $f = 0$.
\end{itemize}
This shows $C_\fg(e) = E_1$.

We finally rescale and replace $h_1$ by $ \frac{1}{2} h_1$.
\end{proofclaim}

\begin{step}\label{p:3:identification:st:coordinatisation}
Coordinatisation.
\end{step}
\begin{proofclaim}
Let $\fc = E_0$, $\fu_+ = E_2$, and $\fu_- = E_{-2}$. By Step~\ref{p:3:identification:st:weights}, there are $e_1 \in \fu_+$ and $f_1 \in \fu_-$ with $[e_1, f_1] = h_1$.

We simultaneously coordinatise $\fc$ and $\fu_+$ as follows. Linearising in dimension $1$, there is a definable field $\bK$ such that $\fc \simeq \bK \Id$ acts by \emph{doubled} scalars on $\fu_+ \simeq \bK_+$, viz.~$[h_\lambda, e_\mu] = 2 e_{\lambda\mu}$.
Notice that $h_1$ was defined consistently.
Since the multiplicative unit is not fixed, we choose it to be $e_1$ and consistently write $\fu_+ = \{e_\lambda: \lambda \in \bK\}$. 
Now we let $f_\lambda = -\frac12 [h_\lambda, f_1]$, so $\fu_- = \{f_\lambda: \lambda \in \bK\}$. Notice that $f_1$ was defined consistently.

We prove $\fg \simeq \fsl_2(\bK)$. It suffices to check a couple of identities. First, by definition, $[h_\lambda, e_\mu] = \frac12 [h_\lambda, [h_\mu, e_1]] = [h_\mu, e_\lambda]$; and $[h_\lambda, f_\mu] = [h_\mu, f_\lambda]$ likewise.

Then $[e_\mu, f_1] = h_\mu$, because:
\[
[[e_\mu, f_1], e_1] = [e_\mu e_1, f_1] + [e_\mu, f_1 e_1] = [h_1, e_\mu] = 2 e_\mu = [h_\mu, e_1],\]
meaning that $[e_\mu, f_1] - h_\mu \in C_\fc(e_1) = 0$ by Step~\ref{p:3:identification:st:weights}.


It follows $[e_\mu, f_\nu] =  h_{\mu\nu}$. Indeed,
\[-2 [e_\mu, f_\nu] = [e_\mu, [h_\nu, f_1]] = [e_\mu h_\nu, f_1] + [h_\nu, e_\mu f_1] = - 2 [e_{\mu\nu}, f_1] = - 2 h_{\mu\nu}.
\]

Last, $[h_\mu, f_\nu] = - 2 f_{\mu\nu}$ since:
\[
[[h_\mu, f_\nu], e_1]  = [h_\mu e_1, f_\nu] + [h_\mu, f_\nu e_1] = 2 [e_\mu, f_\nu] = 2 h_{\mu\nu} = - 2 [f_{\mu\nu}, e_1],
\]
implying $[h_\mu, f_\nu] + 2 f_{\mu\nu}\in C_{\fu_-}(e_1) = 0$ by Step~\ref{p:3:identification:st:weights} again.

The isomorphism $\fg \simeq \fsl_2(\bK)$ is now obvious.
\end{proofclaim}
This completes the recognition and proves Theorem~\ref{t:3}.
\end{proof}

The extension problem is easy, although not considered by Rosengarten.

\begin{corollary}\label{c:3}
Let $\fg$ be a connected, non-soluble Lie ring of Morley rank $3$ and characteristic $\neq 2, 3$. Then $\fg \simeq \fsl_2(\bK)$ for some definable field $\bK$.
\end{corollary}
\begin{proof}
It is enough to show simplicity of $\fg$, which we shall derive from the isomorphism type of $\fg/Z(\fg)$. Let $Z = Z(\fg)$, a definable ideal of $\fg$, and $\overline{\fg} = \fg/Z$.

We first prove that $\overline{\fg}$ is simple. It is enough to prove that it is definable simple (see \S~\ref{s:facts}).
Let $I \triangleleft \fg$ be any proper definable ideal; then $I^\circ \triangleleft\fg$ as well.
If $I^\circ > 0$ then by Theorem~\ref{t:2} both $I^\circ$ and $\fg/I^\circ$ are soluble, against the assumption. Therefore $I$ is finite. By connectedness, $[I, \fg] = 0$ so $I$ is central. This shows that the largest proper definable ideal of $\fg$ is $Z$. In particular, $\overline{\fg}$ is definably simple; being non-abelian, it is simple.

We now conclude $Z = 0$.
By Theorem~\ref{t:3}, $\overline{\fg} \simeq \fsl_2(\bK)$ for some definable field $\bK$.
So there is $\eta \in \overline{\fg}$ such that $\overline{\fg} = E_{-1}(\eta) \oplus E_0(\eta) \oplus E_1(\eta)$. Lift $\eta$ to $h \in \fg$; lifting eigenspaces, one has $\fg = E_{-1}(h) \oplus E_0(h) \oplus E_1(h)$, and each has dimension $1$. It now follows that $Z \leq E_0(h)$. Moreover, $[E_{-1}(h), E_1(h)] = E_0(h)$. Let $z \in Z$. Then there are $a_{-1} \in E_{-1}$ and $a_1 \in E_1$ with $z = [a_{-1}, a_1]$. Let $\alpha_i = (a_i \mod Z) \in E_i(\eta)$. Then $[\alpha_{-1}, \alpha_1] = 0$ in $\fsl_2(\bK)$, so at least one $\alpha_i$ is zero. Thus $a_i \in E_i \cap Z = 0$ and $z = 0$. This proves $Z = 0$ and $\fg \simeq \overline{\fg}$.
\end{proof}

\subsection{Dimension 4}\label{s:4}

\begin{theorem}\label{t:4}
There is no simple Lie ring of Morley rank $4$ and characteristic $\neq 2, 3$.
\end{theorem}

\paragraph{The main idea}
is elementary: use Borel subrings as flintstones. We need at least two, and at least one of them should be large enough. It is unreasonable to have more upper-triangular matrices than lower-triangular ones, and this is the contradictory picture. For this we need two Borel subrings of distinct dimensions intersecting over a $1$-dimensional, `diagonal' subring, creating friction between unbalanced dimensions.

So we shall produce a $3$-dimensional Borel subring $\fb_1$, and another Borel subring $\fb_2$. With rudimentary proxies for `a toral subring' and `the nilpotent radical', $\fb_1$ and $\fb_2$ will intersect in a $1$-dimensional, toral subring $\fc$, and have disjoint nilpotent radicals $I_1 \neq I_2$; also we shall see $\dim I_1 = 2$. Then $[I_1, I_2]$ will be too large to fit into $\fc$. This creates centralisation phenomena between nilpotent radicals in a tight configuration, where nilpotence is not for sharing: a contradiction.

Of course the above strategy has to be implemented without weight theory, since there is no linear structure to begin with. So we are left with basic means inspired by the group-theoretic analysis of $\PGL_2(\bK)$ in terms of Borel subgroups and unipotent radicals \cite{CJTame, DGroupes, DJInvolutive}.
In short the proof does \emph{not} proceed by picking a Cartan subring and studying the root decomposition. To do this one would need a number of tools yet dubious; see final questions in \S~\ref{S:questions}.

\paragraph{The proof starts.}
Let $\fg$ be a simple Lie ring of Morley rank $4$.
There are three main stages.
We first prove minimal simplicity of $\fg$ in Proposition~\ref{p:4:minimalsimple}.
We derive that $\fg$ cannot be a sum of integral weight spaces in Proposition~\ref{p:4:weights}; Proposition~\ref{p:4:E0E1} is an important consequence.
Changing topic, we remove nilpotent Borel subrings in Proposition~\ref{p:4:nonilpotent}; Proposition~\ref{p:4:3dimensionalBorel} then allows us to find a $3$-dimensional Borel subring, whose fine structure is described in Proposition~\ref{p:4:b3}. Proposition~\ref{p:4:contradiction} gives the final contradiction.

\begin{proposition}\label{p:4:minimalsimple}
$\fg$ is minimal simple, viz.~every definable, connected, proper subring is soluble.
\end{proposition}
\begin{proof}
Let $\fh \sqsubset \fg$ be a definable, connected, proper subring.
By Theorem~\ref{t:2}, we may assume $\dim \fh = 3$.
Let $I = C_\fh(\fg/\fh)$, an ideal of $\fh$. Notice $I < \fh$ since otherwise $\fh \triangleleft \fg$, a contradiction. In particular $I^\circ$ is soluble by Theorem~\ref{t:2}. Linearising in dimension $1$ forces abelianity of $\fh/I$. Thus $\fh$ is soluble.
\end{proof}

Much later, Proposition~\ref{p:4:b3} will return to the analysis of $3$-dimensional subrings.

\paragraph{Lack of weights.}

We shall now forbid $\fg$ to be a sum of integral weight spaces.
While mostly trivial, the proof needs attention as we also allow for characteristic $5$. There, dimension $4$ comes extremely close to `Witt-type' root distribution and requires some care. Moreover, the final case of two Borel subrings intersecting over a $2$-dimensional Cartan subring (Step~\ref{p:4:weights:contradiction}) is a key configuration to eliminate.

\begin{proposition}\label{p:4:weights}
Let $h \neq 0$. Then $\fg$ is not a sum of eigenspaces, viz.~$\sum_{i \in \bF_p} E_i(h) < \fg$.
\end{proposition}
\begin{proof}
Unless otherwise indicated, \emph{all eigenspaces in this proof are with respect to $h$.} Let $p$ be the characteristic; by assumption $p \neq 2, 3$, but $5$ is allowed.

Suppose $\sum_{i \in \bF_p} E_i = \fg$.
Bear in mind that the sum is direct.
\begin{step}\label{p:4:weights:-101k}
One has $h\in E_0$. Allowing trivial terms, we may suppose $\fg = E_{-1} \oplus E_0 \oplus E_1 \oplus E_k$ for some $k \notin \{-1, 0, 1\}$.
\end{step}
\begin{proofclaim}
A priori, and up to allowing trivial terms, $\fg = E_i \oplus E_j \oplus E_k \oplus E_\ell$ for $i, j, k, \ell$ \emph{distinct} integers modulo $p$.

We claim $h \in E_0 \neq 0$ and $\fg = E_0 \oplus E_i \oplus E_j \oplus E_k$. Indeed, write $h = e_i + e_j + e_k + e_\ell$ in obvious notation; we may assume $e_\ell \neq 0$. Applying $\ad_h$, we find $ie_i + je_j + ke_k + \ell e_\ell = 0$. In particular, $\ell e_\ell = 0$ and therefore $\ell = 0$. Since the integers were distinct, $e_i = e_j = e_k = 0$ and $h = e_\ell \in E_\ell = E_0$.

We may assume $E_i$, $E_j$, $E_k$ all non-trivial and proper. Properness is since $h \in E_0$. Now if $E_i = 0$, then let $h' = \frac{1}{j} h$. Notice how $E_\ell(h') = E_{j\ell}(h)$, so $\fg = E_0(h') \oplus E_1(h') \oplus E_{j^{-1} k}(h')$. Up to considering $h'$ instead of $h$ (which we call \emph{up to dividing}) we have the desired form.

We may assume $i + j$, $i + k$, $j+k$ all non-zero. Indeed if $j = -i$, then up to dividing we have the desired form.

We may assume $i + j \neq k$ and $i + k \neq j$. Indeed if $i + j - k = i + k - j = 0$, then $2i = 0$ and $i = 0$; a contradiction. So at most one the three equalities of type $i + j = k$ holds.

We may assume $2i \in \{j, k\}$. Indeed since $i + j \neq 0$ and $i + j \neq k$, we have $i + j \notin \{0, i, j, k\}$, so $[E_i, E_j] \leq E_{i + j} = 0$. By our assumptions, $[E_i, E_k] = 0$ likewise. So $0 < E_i < \fg$ is normalised by $E_0 + E_j + E_k$. Since $E_i$ is no ideal of $\fg$, we therefore have $[E_i, E_i] \neq 0$ and $2i \in \{0, i, j, k\}$. But $2i \notin \{0, i\}$ since $i \neq 0$.

We have $2j \in \{i, k\}$ and $2k \in \{i, j\}$. For suppose $j + k = i$. With $2i = j$ this gives $i + k = 0$, a contradiction; the case $2i = k$ is similarly excluded. So $j + k \neq i$, and the previous paragraph also gives $2j \in \{i, k\}$ and $2k \in \{i, j\}$.

Conclusion. We may assume $2i = j$; the other case is similar. If $2j = i$ then $4 = 1$, against $p \neq 3$. So $2j = k$. Similarly, $2k = i$. Therefore $8 = 1$ and $p = 7$. Up to dividing for readability, we now have $\fg = E_0 \oplus E_1 \oplus E_2 \oplus E_4$. But now $E_3 = E_5 = E_6 = 0$ and $E_8 = E_1$, so $E_1 \oplus E_2 \oplus E_4 \triangleleft \fg$, a contradiction.\end{proofclaim}

\begin{step}\label{p:4:weights:-101}
We may assume $\fg = E_{-1} \oplus E_0 \oplus E_1$. Moreover $E_{-1}$, $E_0$, $E_1$ are non-trivial and proper.
\end{step}
\begin{proofclaim}
By Step~\ref{p:4:weights:-101k} we may assume $\fg = E_{-1} \oplus E_0 \oplus E_1 \oplus E_k$ for some $k \notin \{-1, 0, 1\}$.

We first show that we may assume $k = 2$.
If $k = -2$, then up to considering $-h$, we have the desired form. Therefore suppose $k \neq \pm 2$, so that $k \notin \{0, \pm 1, \pm 2\}$. In particular $E_{-2} = E_2 = 0$, implying that $E_{-1}$ and $E_1$ are abelian subrings. Let $\fh = E_{-1} + E_0 + E_1$, a definable, connected subring of $\fg$.
If $\fh = \fg$, then up to taking $k = 2$ as a dummy we have the desired form. So we may assume $\fh < \fg$. By minimal simplicity (Proposition~\ref{p:4:minimalsimple}), $\fh$ is soluble.
By Step~\ref{p:4:weights:-101k}, $h \in E_0 \leq \fh$, so $E_1 = [h, E_1] \leq \fh'$, and $E_{-1} \leq \fh'$ likewise. By solubility, $\fh' < \fh$, implying $[E_{-1}, E_1] < E_0$. If $\dim E_0 = 2$ then since $\fh < \fg$ is proper, at least one of $E_{-1}, E_1$ is trivial; hence $[E_{-1}, E_1] = 0$. If $\dim E_0 = 1$ then $[E_{-1}, E_1] = 0$. So $[E_{-1}, E_1] = 0$ in either case. Moreover $k \pm 1 \notin \{-1, 0, 1, k\}$, so $E_{-1}$ and $E_1$ are ideals of $\fg$. Since $E_0 \neq 0$ they are proper, hence trivial. There remains $\fg = E_0 \oplus E_k$ and $E_k$ is a proper ideal of $\fg$, hence trivial. Finally $\fg = E_0$, so $h \in Z(\fg) = 0$, a contradiction. So we may assume $k = 2$.
Hence $\fg = E_{-1} \oplus E_0 \oplus E_1 \oplus E_2$.

We prove that $E_{-1}, E_0, E_1$ are non-trivial and proper.
\begin{itemize}
\item
We already know $0 < E_0 < \fg$.
\item
If $E_{-1} = 0$ then (even when $p = 5$) $E_2$ is an ideal of $\fg$; hence trivial. Then $E_1$ is an ideal of $\fg$; hence trivial as well, and this is a contradiction.
\item
If $E_1 = 0$ then $[E_{-1}, E_2] = 0$ so (even when $p = 5$) $E_{-1}$ is an ideal of $\fg$; hence trivial, a contradiction.
\end{itemize}

It remains to show that $E_2 = 0$. Suppose not. Then $\fg = E_{-1} \oplus E_0 \oplus E_1$, and each has dimension $1$.

We claim $[E_{-1}, E_1] = E_0$, $[E_1, E_1] = E_2$, and $[E_{-1}, E_2] = E_1$.
\begin{itemize}
\item
If $[E_{-1}, E_1] = 0$ then (even when $p = 5$) $E_{-1} + E_1 + E_2$ is an ideal of $\fg$, a contradiction.
\item
If $[E_1, E_1] = 0$, then $\fh = E_{-1} + E_0 + E_1$ is a proper subring, hence soluble by Proposition~\ref{p:4:minimalsimple}. But $h \in E_0$ so $\fh' \geq E_{-1} + E_1 + [E_{-1}, E_1] = \fh$, a contradiction.
\item
Suppose $[E_{-1}, E_2] = 0$. Write $h \in E_0 = [E_{-1}, E_1]$ as $h = [e_{-1}, e_1]$; let $e_2 \in E_2$. Then:
\[2 e_2 = [h, e_2] = [[e_{-1}, e_1], e_2] = [e_{-1} e_2, e_1] + [e_{-1}, e_1 e_2] = 0,\]
a contradiction. (If the characteristic is not $5$ there is a shorter argument: $E_2$ is an ideal of $\fg$, a contradiction.)
\end{itemize}

We claim that for any $x_0 \in E_0$ there is $y_2 \in E_2\setminus\{0\}$ such that $[x_0, y_2] = 0$.
Let $e_{-1} \in E_{-1}$ be such that $[e_{-1}, E_1] = E_0$ and $[e_{-1}, E_2] = E_{1}$; there is one such since otherwise $C_{E_{-1}}(E_1) \cup C_{E_{-1}}(E_2) = E_{-1}$, a contradiction. Let $x_0 \in E_0$. Then there is $y_1 \in E_1$ with $[y_1, e_{-1}] = x_0$. Now there is $y_2 \in E_2$ with $[e_{-1}, y_2] = y_1$. On the other hand, $[E_1, E_2] = 0$. Altogether, one has:
\begin{align*}
[x_0, y_2] & = [[y_1,e_{-1}],y_2]\\
& = [y_1 y_2, e_{-1}] + [y_1, e_{-1} y_2]\\
& = [0, e_{-1}] + [y_1, y_1]\\
& = 0.
\end{align*}

Conclusion. Since $h \in E_0$, the action of $E_0$ on $E_2$ is non-trivial. Linearising in dimension $1$, there is a finite ideal $F \triangleleft E_0$ such that $E_0/F$ acts by scalars, hence freely, on $E_2$. This contradicts the last claim.
\end{proofclaim}

\begin{step}\label{p:4:weights:contradiction}
Contradiction.
\end{step}
\begin{proofclaim}
By Step~\ref{p:4:weights:-101}, exactly one of $E_{-1}, E_0, E_1$ has dimension $2$; moreover $E_{\pm 1}$ are abelian.
For the end of the argument, bear in mind that $E_{-1}$ does not normalise $E_1$ nor vice-versa, since the normalised one would be an ideal of $\fg$.
There are essentially two cases.
\begin{itemize}
\item
Suppose $\dim E_1 = 2$.
Let $a \in E_{-1}$. Since $\ad_a$ takes $E_1$ to $E_0$, the subgroup $X_a = C_{E_1}^\circ(a)$ is infinite; since $E_1$ is abelian, $X_a$ is actually a subring. Notice that both $a$ and $h$ normalise $X_a$. Moreover $C_\fg^\circ(X_a) \geq E_1$. If $C_\fg^\circ(X_a) > E_1$, then it is a $3$-dimensional Borel subring. In that case, by Lemma~\ref{l:selfnormalisation}, $h \in C_\fg^\circ(X_a)$, a contradiction as the action of $h$ on $X_a$ is non-trivial. Hence $C_\fg^\circ(X_a) = E_1$ is normalised by $a$. This shows that $E_{-1}$ normalises $E_1$, a contradiction.
%

This rules out $\dim E_{-1} = 2$ as well, by considering $-h$.
\item
Suppose $\dim E_0 = 2$. (It is unclear whether $E_0$ is abelian.)
Since $E_0$ contains $h$, it is trivial neither on $E_1$ nor on $E_{-1}$. Linearising in dimension $1$, both $\ft_1 = C_{E_0}^\circ(E_1)$ and $\ft_{-1} = C_{E_0}^\circ(E_{-1})$ are subrings of dimension $1$. In particular they are abelian.

We show $\ft_1 \neq \ft_{-1}$. Otherwise denote it by $\ft_{\pm 1}$ and let $\fh = C_\fg^\circ(\ft_{\pm 1})$, a proper subring of $\fg$. By assumption, $\fh \geq E_{-1} + E_1 + \ft_{\pm 1}$; so equality holds. In particular $[E_{-1}, E_1] \leq (\fh\cap E_0)^\circ = \ft_{\pm 1}$. Now $h$ normalises $E_{-1}$, $E_1$, and their bracket $\ft_{\pm 1}$; so $h$ normalises $\fh$. By Lemma~\ref{l:selfnormalisation}, one gets $h \in \fh = E_{-1} \oplus E_0^\fh \oplus E_1$, and finally $h \in E_0^\fh = \ft_{\pm 1}$, a contradiction.

Let $\fc = [E_{-1}, E_1] \neq 0$. Since $E_0$ normalises $\fc$, this is a subring.
We prove $\ft_1 \not\leq \fc$; suppose inclusion holds. Then $\ft_1 \neq \ft_{-1}$ acts non-trivially on $E_{-1}$, so linearising in dimension $1$, there is $t_1 \in \ft_1$ acting on $E_{-1}$ as $1$; by definition it acts trivially on $E_1$. Then by assumption, $t_1 \in \fc = [E_{-1}, E_1]$. At most two elements $e_{-1}, e_{-1}' \in E_{-1}$ suffice to write $\fc = \ad_{e_{-1}}(E_1) + \ad_{e_{-1}'}(E_1)$. So there are $e_1, e_1' \in E_1$ with $t_1 = [e_{-1}, e_1] + [e_{-1}', e_1']$. But then:
\begin{align*}
0 & = [t_1, t_1] = [t_1, [e_{-1}, e_1]] + [t_1, [e_{-1}', e_1']]\\
& = [t_1 e_{-1}, e_1] + [e_{-1}, t_1 e_1] + [t_1 e_{-1}', e_1'] + [e_{-1}', t_1 e_1']\\
& = [e_{-1}, e_1] + [e_{-1}', e_1'] = t_1,
\end{align*}
a contradiction. Notice $\ft_{-1} \not\leq \fc$ for the same reasons.

Conclusion. This implies $\dim \fc = 1$. Indeed, $\dim \fc = 2$ would contradict $\ft_1 \not\leq \fc$.
Finally let $\fh = E_{-1} + \fc + E_1$, a subring. If $\fc$ centralises $E_1$, then $\fc = \ft_1$, a contradiction. For the same reason, $\fc$ does not centralise $E_{-1}$. So $\fh' = \fh$, against Proposition~\ref{p:4:minimalsimple}.
\qedhere
\end{itemize}
\end{proofclaim}

This proves Proposition~\ref{p:4:weights}.
\end{proof}

As a consequence we derive a severe control on sizes of weight spaces.

\begin{proposition}\label{p:4:E0E1}
Suppose $h \in \fg$ is such that $E_0(h) \neq 0$ and $E_1(h) \neq 0$. Then $\dim E_0(h) = \dim E_1(h) = 1$ and $E_2(h) = 0$.
\end{proposition}
\begin{proof}
\emph{All eigenspaces in this proof are with respect to $h$.}

Let $\fb = E_0 + E_1 + E_2$, a priori a definable, connected subgroup. By Proposition~\ref{p:4:weights}, $E_3 = E_4 = 0$, so $\fb$ is actually a proper subring. We must prove $\dim \fb = 2$ and suppose not; then $\fb$ is a $3$-dimensional Borel subring. By Lemma~\ref{l:selfnormalisation}, it is self-normalising in $\fg$. Therefore $h \in N_\fg(\fb) = \fb$.

Since $h \in \fb$, for $i \in \{1, 2\}$ one has $E_i  = [h, E_i] \leq \fb'$. It follows $\fb' = E_1 + E_2 + E_0^{\fb'}$, where each of the last two terms is allowed to be trivial. Then $[h, \fb'] \leq E_1 + E_2$.

Consider the action of $\fb$ on $\fg/\fb$. Linearising in dimension $1$, $\fb'$ centralises $\fg/\fb$. So for $i \in \{1, 2\}$ and $e \in E_i$ one has $B_e \leq \fb$. By Proposition~\ref{p:4:weights} again, $E_{-i} = 0$, so $\ad_h + i$ is surjective and:
\[[h, B_e] = \im \ad_h \ad_e = \im \ad_e (\ad_h + i) = \im \ad_e = B_e.\]
Since $B_e \leq \fb$ and $h \in \fb$ one has $B_e = [h, B_e] \leq \fb'$. Now $B_e = [h, B_e] \leq [h, \fb'] \leq E_1 + E_2$. Summing over $e \in E_i$ and $i \in \{1, 2\}$, for $a \in E_1 + E_2$ one has $B_a \leq E_1 + E_2$. Hence $E_1 + E_2 \triangleleft \fg$. But $E_0 \neq 0$ and $E_1 \neq 0$: a contradiction to simplicity.
\end{proof}

\paragraph{A first analysis of Borel subrings.}
We now change thread and study abstract Borel subrings.

\begin{proposition}[see Proposition~\ref{p:3:nonilpotent}]\label{p:4:nonilpotent}
No Borel subring is nilpotent.
\end{proposition}
\begin{remark}
Some of the ideas in the following proof could be pushed further. One suspects that if $\fg$ is minimal simple, then no Borel subring has dimension $1$.
\end{remark}
\begin{proof}
Towards a contradiction, let $\fb$ be one such. Let $Z = Z^\circ(\fb) > 0$. It is a possibility that $Z = \fb$. For $z \in Z\setminus\{0\}$, one has $\fb \sqsubseteq C_z$; by minimal simplicity, equality holds.

\begin{step}[see Step~\ref{p:3:nonilpotent:st:1modules}]\label{p:4:nonilpotent:st:1modules}
Let $X$ be a $1$-dimensional subquotient $\fb$-module. Then $\fb$ centralises $X$. In particular, $\dim \fb \leq 2$.
\end{step}
\begin{proofclaim}
The former claim implies the latter. Indeed, if $\dim \fb = 3$ then $\fb$ acts trivially on $\fg/\fb$, so $[\fb, \fg] \leq \fb$ and $\fb$ is an ideal of $\fg$, a contradiction. So we deal with the former claim.

We first prove that if $Y$ is a $1$-dimensional subquotient $Z$-module, then $Z$ centralises $Y$. (Both the conclusion \emph{and assumption} are weaker.) Suppose not. Linearising in dimension $1$, there is $z \in Z$ acting on $Y$ like $\Id_Y$. \emph{In the present paragraph, eigenspaces are with respect to $z$.}
Lifting the eigenspace, $E_1 \leq \fg$ is non-trivial. Since $p \geq 5$, the multiplicative order of $2 \in \bF_p^\times$ is at least $3$; on the other hand $E_0 > 0$. So there is $\ell \in \bF_p^\times$ such that $E_\ell \neq 0$ but $E_{2\ell} = 0$. Then $E_0 + E_\ell$ is a subring properly containing $E_0 = C_z = \fb$. It follows $\fg = E_0 + E_\ell$ and $E_\ell \triangleleft \fg$, a contradiction. So the claim about $1$-dimensional $Z$-subquotient modules holds. (We lost $z$.)

To prove our full claim, we may suppose that $\fb$ is non-abelian, implying $\dim \fb \geq 2$, and that $\fb$ does not act trivially on $X$. Linearising in dimension $1$, there is $h \in \fb$ acting on $X$ as the identity. \emph{From now on, all eigenspaces are with respect to $h$.}
Lifting the eigenspace, there is $h \in \fb$ with $E_1 \leq \fg$ non-trivial. We aim for a contradiction.

We claim that $\fg = \fb \oplus E_1$ where $\dim \fb = \dim E_1 = 2$, and $\fg/\fb \simeq E_1$ as $Z$-modules. (Not `as $\fb$-modules', since $\fb$ need not normalise $E_1$.)
First, $\fb \cap E_1 = 0$. This is because $\fb \cap E_1$ is contained in each term of the descending nilpotence series, which vanishes by nilpotence of $\fb$.
Second, $\dim E_1 \geq 2$. Otherwise $\dim E_1 = 1$ and the action of $Z$ on $E_1$ is trivial, so $E_1 \leq C_\fg^\circ(Z) = \fb$, a contradiction. In particular, $\dim \fb = \dim E_1 = 2$.
Now $\fg = \fb \oplus E_1$ and $E_1 \simeq \fg/\fb$ as $Z$-modules.

We show that $Z$ acts freely and irreducibly on $\fg/\fb$.
If $E_1$ is not $Z$-irreducible, there is a $1$-dimensional $Z$-module $Y < E_1$. Now $Z$ acts trivially on $Y$ so $Y \leq C_\fg^\circ(Z) = \fb$ and $Y \leq E_1 \cap \fb$, a contradiction again. Hence $E_1$ is $Z$-irreducible. If some $z \in Z\setminus\{0\}$ acts trivially then $E_1 \leq C_z = \fb$, a contradiction. Linearising the abelian action, there is a field $\bK$ with $E_1 \simeq \bK_+$ and $Z \hookrightarrow \bK \Id_{E_1}$.
Therefore $Z$ acts freely on $\fg/\fb \simeq E_1$.

We finish the proof by turning to the action of $\fb$ on $\fg/\fb$.
Let $I = C_\fb(\fg/\fb)$. Since $Z$ acts freely on $\fg/\fb$, one has $I \cap Z = 0$. In particular $\dim I < 2$. If $\dim I = 1$ then $\fb = I + Z$ is a sum of two abelian ideals, hence abelian itself: a contradiction.
Therefore $\dim I = 0$.
Notice that $\fg/\fb$ is $\fb$-irreducible as it was $Z$-irreducible.
Now $ZI/I$ is infinite, and can be used to linearise the action of $\fb/I$ (see \S~\ref{s:facts}). Hence $\fb/I$ is a linear nilpotent Lie ring, and $\fg/\fb$ is $\fb$-irreducible. This forces the linear dimension to be $1$, and therefore $\fb/I$ is abelian. But then $\fb'$ is finite and connected: $\fb$ is abelian, a contradiction.
\end{proofclaim}

Let $W = \fg/\fb$, a $\fb$-module, and $I = C_\fb(W)$, an ideal of $\fb$. The proof of Proposition~\ref{p:3:nonilpotent} cannot be followed literally but some ideas will look familiar.

\begin{step}
$\dim \fb = 1$; in particular, $\fb$ is abelian.
\end{step}
\begin{proofclaim}
Suppose $\dim \fb = 2$.

We claim that $W$ is $\fb$-irreducible. Otherwise there is a $3$-dimensional $\fb$-module $V_3$ with $\fb < V_3 < \fg$. By Step~\ref{p:4:nonilpotent:st:1modules}, $\fb$ acts trivially on $V_3/\fb$, so $V_3 \leq N_\fg^\circ(\fb)$: a contradiction.

Now we show $\dim I = 1$ and $I^\circ \leq Z$. The former implies the latter. Indeed if $\fb$ is abelian, we are done; while if $\fb$ is not, then its only $1$-dimensional, connected ideal is $Z$, so $I^\circ = Z$. So we focus on proving $\dim I = 1$.
If $\dim I = 2$ then $\fb$ is an ideal of $\fg$: a contradiction.
We suppose $\dim I = 0$ and give a contradiction. First, $\fb/I$ can be linearised thanks to $ZI/I > 0$. Since $\fb$ is nilpotent and $W$ is $\fb$-irreducible, this forces the linear dimension to be $1$ and $\fb$ to be abelian. Moreover there is a field structure with $\fb/I\simeq \bK\Id_W$ and $W \simeq \bK_+$. In particular there is $h \in \fb$ acting on $\fg/\fb$ as the identity. We lift the eigenspace; thus $E_1 > 0$, but also $E_0 = \fb$ by abelianity. We prove $\dim E_1 \geq 2$ as follows (covering properties are unclear to us; $\dim E_1^{\fg/\fb}(h) = 2$ does not seem to be a sufficient argument). Let $z \in Z\setminus \{0\}$; for $x \in \fg$ there is $b \in \fb$ such that $[h, x] = x + b$, so:
\[[h, [z, x]] = [hz, x] + [z, hx] = [z, x + b] = [z, x].\]
Therefore $\ad_z\colon \fg \to E_1$. However $\ker^\circ \ad_z = C_z = \fb$, so $\dim E_1 \geq \dim \fg - \dim \fb = 2$, as wished.
Equality follows as $\fb \cap E_1 = 0$ by abelianity. So $\fg = \fb \oplus E_1$ and $\fb = E_0$ normalises $E_1$, an ideal of $\fg$; this is a contradiction. The claim is proved.

It follows $N_\fg(\fb) = \fb$. Indeed, $\fb/I$ is abelian so we may linearise the action of $\fb/I$ on $\fg/\fb$: it is by scalars, hence free. Now if $n \in N_\fg(\fb)$ and $b \in \fb \setminus I$, then $[b, n] \in \fb$ implies $n \in \fb$. (This is the argument of Lemma~\ref{l:selfnormalisation}, which we cannot invoke as stated.)

We derive a contradiction.
Let $z \in I^\circ \leq Z$ be non-trivial. Since $C_z = \fb$, one has $\dim B_z = 2$. By definition of $I$, one also has $B_z \leq \fb$, so equality follows.
Hence there is $x \in \fg$ with $z = [x, z]$. For arbitrary $b \in \fb$ one has:
\[[[b, x], z] = [bz, x] + [b, xz] = [b, z] = 0,\]
so $[\fb, x] \leq C_z = \fb$. Hence $x \in N_\fg(\fb) = \fb$, and $z = [x, z] = 0$, a contradiction.
\end{proofclaim}

\begin{step}[see Step~\ref{p:3:nonilpotent:st:contradiction}]\label{p:4:nonilpotent:st:contradiction}
Contradiction.
\end{step}
\begin{proofclaim}
Throughout, $h$ will stand for an arbitrary non-zero element of $\fb$. Let $B_\fb = [\fb, \fg]$, a $\fb$-module containing each $B_h$ for $h \in \fb\setminus\{0\}$.

We first show that $B_h = B_\fb$ and $\dim B_h = 3$. Indeed, by abelianity, $B_h$ is $\fb$-invariant; now $C_h = \fb$, so $\dim (\fg/B_h) = 1$. By Step~\ref{p:4:nonilpotent:st:1modules}, $\fb$ acts trivially on $\fg/B_h$. So $[\fb, \fg] \leq B_h$. Hence for any $h' \in \fb$ one has $B_{h'} \leq B_h$ and equality follows. Thus $B_h$ does not depend on $h \in \fb\setminus\{0\}$; it equals $B_\fb$.

We claim that $C_h^2 = C_h = \fb$ and $B_h^2 = B_\fb$. Indeed $C_h = \fb$ is abelian, so by Lemma~\ref{l:Cn}, $C_h^2$ is a subring of $\fg$. But $\dim C_h^2 \leq 2 \dim C_h < \dim \fg$, so $C_h^2$ is a proper subring containing $C_h = \fb$. Hence $C_h^2 = \fb = C_h$ and therefore $B_h^2 = B_h = \fb$. 

We deduce that $\fg = \fb + B_\fb$ and the action of $\fb$ on $B_\fb$ is faithful. Otherwise $\fb \leq B_\fb$; hence $C_h \leq B_h$, which forces $B_h^2 < B_h$, a contradiction.
Moreover, if $h \in \fb$ centralises $B_\fb$, it also centralises $\fb + B_\fb = \fg$ and is therefore zero.

We now contend that $B_\fb$ is $\fb$-irreducible. Suppose it is not. If there is a $1$-dimensional $\fb$-submodule $V_1 < B_\fb$, then by Step~\ref{p:4:nonilpotent:st:1modules} one has $V_1 \leq C_\fg^\circ(\fb) = \fb$, so $\fb = V_1 \leq B_\fb$, a contradiction. If there is a $2$-dimensional $\fb$-submodule $V_2 < B_\fb$, then by Step~\ref{p:4:nonilpotent:st:1modules} again, one has $[\fb, B_\fb] \leq V_2$, so $B_\fb^2 < B_\fb$, a contradiction.

We may thus linearise the action of $\fb$ on $B_\fb$; there is a field $\bK$ such that $B_\fb \simeq \bK_+$ and $\fb$ acts by scalars on $B_\fb$. Let us fix some notation: there are $\varphi\colon B_\fb \simeq \bK_+$ and $\chi\colon \fb \hookrightarrow \bK_+$ such that, for $h \in \fb$ and $x \in B_\fb$:
\[\varphi([h, x]) = \chi(h) \cdot \varphi(x).\]

Fix arbitrary $x \in B_\fb$ and let $A_n = [\fb, \dots, [\fb, x]{\scriptscriptstyle \dots}]$.
Also let $A'_n = \varphi(A_n)$.
These definable, connected subgroups of $\bK_+$ need not form an ascending series but their dimensions do. Let $n$ be such that $\dim A'_n = \dim A'_{n+1}$.
If $h \in \fb$ and $q \in A'_n$, say $q = \varphi(y)$ with $y \in A_n$, then:
\[\chi(h) \cdot q = \varphi([h, y]) \in A'_{n+1}.\]
But dimensions match, so $\chi(h) A'_n = A'_{n+1}$ as subgroups of $\bK_+$.

So fix $h_0 \in \fb\setminus\{0\}$. For arbitrary $h_1 \in \fb$, one has $\frac{1}{\chi(h_0)} \chi(h_1) \in N_\bK(A'_n) = \{\lambda \in \bK: \lambda A'_n \leq A'_n\}$. The latter is therefore an infinite definable subring, implying $N_\bK(A'_n) = \bK$. Hence $A'_n$ is a non-trivial ideal of a field, and $A'_n = \bK$. Going back through $\varphi$, we have $A_n = B_\fb$. But $[A_n, A_n] \leq B_\fb$ by Lemma~\ref{l:Rosengarten}. So $B_\fb$ is a subring. As it is normalised by $\fb$ and $\fg = \fb + B_\fb$, we have a final contradiction.
\end{proofclaim}

This completes the proof of Proposition~\ref{p:4:nonilpotent}.
\end{proof}

Using non-nilpotence of Borel subrings (Proposition~\ref{p:4:nonilpotent}) and the impossibility of a global weight decomposition (Proposition~\ref{p:4:weights}), we derive the existence of a $3$-dimensional Borel subring.

\begin{proposition}\label{p:4:3dimensionalBorel}
There is a $3$-dimensional Borel subring.
\end{proposition}
\begin{proof}
Suppose not.
The proof has something common with the weight decomposition of Proposition~\ref{p:3:identification}, Step~\ref{p:3:identification:st:weights}.

Let $\fb$ be a Borel subring; by Proposition~\ref{p:4:nonilpotent}, $\dim \fb = 2$ and $\fb$ is non-nilpotent. Using Theorem~\ref{t:2} we let $h \in \fb$ act on $\fb'$ as $\Id_{\fb'}$. Let $e \in \fb'\setminus\{0\}$.

We claim that $C_e = \fb'$. Otherwise $\fb' \sqsubset C_e \sqsubset \fg$, so $C_e$ is a $2$-dimensional Borel subring; it is non-nilpotent by Proposition~\ref{p:4:nonilpotent}. Now $C_e \neq \fb = N_\fg^\circ(\fb')$ so $\fb' \neq C_e'$. By Theorem~\ref{t:2}, $\fb'$ is a Cartan subring of $C_e$. Moreover $h$ normalises $C_e$ since for $c \in C_e$ one has:
\[[[h, c], e] = [he, c] + [h, ce] = [e, c] = 0.\]
Hence $h$ also normalises $C_e'$. Since $\fb'$ is a Cartan subring of the non-nilpotent, $2$-dimensional ring $C_e$, Theorem~\ref{t:2} gives $\eta \in \fb'$ acting on $C_e'$ as $\Id_{C_e'}$. Then for any $a \in C_e'$ one has $[h, a] \in C_e'$ and therefore:
\[
[h, a] = [h, [\eta, a]] = [h \eta, a] + [\eta, ha] = [\eta, a] + [h, a] = a + [h, a].\]
So $a = 0$, a contradiction.

We derive $B_e^2 = \fb$. Indeed $C_e = \fb'$ is abelian. By Lemma~\ref{l:Cn}, $C_e^2$ is a subring of $\fg$; since $\dim C_e^2 \leq 2 \dim C_e = 2$, it is proper. Now $\fb \sqsubseteq C_e^2$ so equality holds. Hence $\fb' = C_e \sqsubset C_e^2 = \fb$. In particular, $\dim B_e^2 = 2$.
We turn to the action of $\fb$ on $\fg/\fb$. If $\fg/\fb$ is $\fb$-irreducible, then by Corollary~\ref{c:Lie2}, $\fb'$ centralises $\fg/\fb$; in particular $B_e^2 \leq B_e \leq \fb$. If $\fg/\fb$ is $\fb$-reducible, we still get $B_e^2 \leq \fb$ by linearising in dimension $1$ \emph{twice}. In either case, $B_e^2 = \fb$.

We now consider eigenspaces with respect to $h$ and show $\fg = E_{-2} + E_{-1} + E_0 + E_1$.
Let us turn to the adjoint action of $h$. First, $E_0 \neq 0$ since $h \in \fb$ which is soluble. By construction, $\fb' \leq E_1$. Suppose $E_{-1} = 0$. Then by surjectivity of $\ad_h +1$ one gets:
\[B_e = \im \ad_e = \im \ad_e (\ad_h + 1) = \im \ad_h \ad_e = [h, B_e].\]
However $B_e \geq B_e^2 = \fb$ intersects $C_h$ by solubility, a contradiction.
We apply just the same argument to prove $E_{-2} \neq 0$; otherwise:
\[B_e^2 = \im \ad_e^2 (\ad_h + 2) = \im \ad_h \ad_e^2 = [h, B_e^2],\]
but $B_e^2 = \fb$, a contradiction.

Thus $\fg = E_{-2} + E_{-1} + E_0 + E_1$, against Proposition~\ref{p:4:weights}.
\end{proof}

\paragraph{The analysis of $3$-dimensional Borel subrings.}
We now combine the two main threads (weights, and large Borel subrings) into the analysis of $3$-dimensional Borel subrings.

\begin{proposition}\label{p:4:b3}
Let $\fb \sqsubset \fg$ be a $3$-dimensional Borel subring. Then:
\begin{enumerate}[label=(\roman*)]
\item\label{p:4:b3:st:I}
$I = C_\fb^\circ(\fg/\fb) \triangleleft \fb$ is a $2$-dimensional, nilpotent ideal of $\fb$ containing $\fb'$;
\item\label{p:4:b3:st:J}
if $J \leq I$ is a $1$-dimensional subring, then $J \leq Z(I)$;
\item\label{p:4:b3:k-l}
if $h \in \fb$ is such that $I$ is a sum of eigenspaces, then either $I \leq E_0$ or $I = E_k \oplus E_{-k}$ is a sum of $1$-dimensional subrings, with $k \neq 0$;
\item\label{p:4:b3:st:Z}
$Z^\circ(\fb) = 0$;
\item\label{p:4:b3:st:buniqueonI}
$\fb$ is the only Borel subring containing $I$;
\item\label{p:4:b3:st:Iabelian}
$I$ is abelian; for $i \in I\setminus Z(\fb)$ one has $C_i = I$.
\end{enumerate}
\end{proposition}
Be however careful that until~\ref{p:4:b3:st:Iabelian} is proved, a subgroup of $I$ need not be a subring.
\begin{proof}\leavevmode
\begin{enumerate}[label=(\roman*)]
\item
By Lemma~\ref{l:preHrushovski}, $I$ is nilpotent; it is proper in $\fb$ since $\fb$ is no ideal of $\fg$. Consider the action of $\fb$ on $\fg/\fb$. Linearising in dimension $1$, $\dim(\fb/I) \leq 1$ and equality holds.
\item
By nilpotence, $Z^\circ(I) \neq 0$.
If $J \not \leq Z^\circ(I)$ then $I = J + Z^\circ(I)$ and $I$ is abelian, a contradiction.
\item
\emph{All eigenspaces in this proof are with respect to $h$.}
Suppose $I \not\leq E_0$, so there is $k \neq 0$ with $E_k^I \neq 0$. Since $E_0^\fb \neq 0$ by solubility, Proposition~\ref{p:4:E0E1} implies $\dim E_0 = \dim E_k = 1$. Therefore write $I = E_k \oplus E_\ell$; possibly $\ell = 0$, but even so $\dim E_\ell = 1$.

We may assume $E_{-k} \neq 0$. Indeed, suppose $E_{-k} = 0$. Fix $e \in E_k = E_k^\fb \leq I$. Notice $B_e \leq \fb$. Now $E_{-k} = 0$ implies $[h, B_e] = \im \ad_h \ad_e = \im \ad_e (\ad_h + k) = \im \ad_e = B_e$. Therefore $B_e \leq \fb' \leq I$. If $\ell = 0$ then $B_e \leq [h, I] = E_k$, and summing over $e \in E_k$ we obtain $E_k \triangleleft \fg$, a contradiction. Hence $\ell \neq 0$. If $E_{-\ell} = 0$, then for $e \in E_\ell$ we also have $B_e \leq I$; summing over $e \in E_k \cup E_\ell$ we obtain $I \triangleleft \fg$, a contradiction. Hence $\ell\neq 0$ and $E_{-\ell} \neq 0$. Up to exchanging, we may assume $E_{-k} \neq 0$.

Let $\fh = E_{-k} + E_0 + E_k$. By Proposition~\ref{p:4:E0E1}, each term has dimension $1$. By Proposition~\ref{p:4:weights}, $E_{\pm 2k} = 0$, so $\fh$ is a $3$-dimensional subring, therefore self-normalising by Lemma~\ref{l:selfnormalisation}. In particular, $h \in \fh$; hence $E_{-k} + E_k \leq \fh'$. But $\fh$ is soluble by Proposition~\ref{p:4:minimalsimple}, so $[E_{-k}, E_k] = 0$. However $I \leq C_\fg^\circ(E_k)$ by~\ref{p:4:b3:st:J}. If $I < C_\fg^\circ(E_k)$, then the latter is a $3$-dimensional subring, hence self-normalising; then $h \in C_\fg^\circ(E_k)$, a contradiction. So $I = C_\fg^\circ(C_k) \geq E_{-k}$, and the claim is proved.
\item
Let $Z = Z^\circ(\fb)$ and suppose $Z \neq 0$.
Since $\fb$ is not nilpotent by Proposition~\ref{p:4:nonilpotent}, $\dim Z = 1$ and $\fb/Z$ is non-nilpotent. In particular $Z \leq I$.
%
By Theorem~\ref{t:2}, there is $\eta \in \fb/Z$ acting on $(\fb/Z)'$ like the identity. We lift to $h \in \fb$ with the same action, and now lift the eigenspace.
\emph{All eigenspaces in this proof are with respect to $h$.} Thus $E_1^\fb \neq 0$.
By Proposition~\ref{p:4:E0E1}, one has $\dim E_0 = \dim E_1 = 1$, so $E_0 = E_0^\fb = Z \leq I$; of course $E_1 = E_1^\fb \leq I$ as well, contradicting~\ref{p:4:b3:k-l}.
%
%
\item
Let $\fb_2 \neq \fb$ be another Borel subring containing $I$. Clearly $\dim \fb_2 = 3$, so our earlier analysis applies. Let $I_2 = C_{\fb_2}^\circ(\fg/\fb_2)$. Since $I_2 \triangleleft \fb_2$ one has $I_2 \neq I$. So $X = (I \cap I_2)^\circ$ has dimension exactly $1$. Now $X$ is a subring of both $I$ and $I_2$, so by~\ref{p:4:b3:st:J} one has $C_\fg^\circ(X) \geq I + I_2$ and therefore $C_\fg^\circ(X)$ is a $3$-dimensional Borel subring with an infinite centre, against~\ref{p:4:b3:st:Z}.
%
\item
Suppose not. Then $J = I' \triangleleft \fb$ is a $1$-dimensional ideal of $\fb$.


There is $h \in \fb$ with $J = E_1(h)$. Indeed, by~\ref{p:4:b3:st:Z}, the action of $\fb$ on $J$ is non-trivial. Linearising in dimension $1$, there is $h \in \fb$ acting on $J$ like the identity. \emph{All eigenspaces are now with respect to $h$.}

We may assume $E_{-1} \neq 0$. If not then for $j \in J = E_1$ one has $[h, B_j] = B_j$ and quickly $B_j \leq I$. Since $J \leq Z(I)$ by~\ref{p:4:b3:st:J}, one also has $C_j \geq I$. Now by~\ref{p:4:b3:st:buniqueonI}, $\fb$ is the only Borel subring containing $I$ and $j \notin Z(\fb)$ because of the action of $h$, so $C_j = I$. It follows $B_j = I = C_j$. By Lemma~\ref{l:divisors}, $B_j = I$ is abelian and we are done.

Let $\fh = E_{-1} + E_0 + E_1$, which by Proposition~\ref{p:4:weights} is a $3$-dimensional Borel subring. By Lemma~\ref{l:selfnormalisation}, $h \in \fh$. In particular $E_{-1} + E_1 \leq \fh' < \fh$, implying $[E_{-1}, E_1] = 0$. However $I \leq C_\fg^\circ(E_1)$ by~\ref{p:4:b3:st:J}, and equality must hold since otherwise $h \in C_\fg^\circ(E_1)$ by Lemma~\ref{l:selfnormalisation}. So $E_{-1} \leq C_\fg^\circ(E_1) = I$ and $I = E_{-1} + E_1$, hence abelian.

It remains to take $i \in I \setminus Z(\fb)$ and prove $C_i = I$. One inclusion is clear. If $C_i > I$ then by~\ref{p:4:b3:st:buniqueonI}, $C_i = \fb$ and $i \in Z(\fb)$, a contradiction.
\qedhere
\end{enumerate}
\end{proof}

\paragraph{The final contradiction.}

\begin{proposition}\label{p:4:contradiction}
The configuration is inconsistent.
\end{proposition}
\begin{proof}
Let $\fb \sqsubset \fg$ be a $3$-dimensional Borel subring, given by Proposition~\ref{p:4:3dimensionalBorel}.
Since $\fb$ is no ideal of $\fg$, the action of $\fb$ on $\fg/\fb$ is non-trivial. Linearising in dimension $1$, there is $h \in \fb$ acting on $\fg/\fb$ like the identity.
\emph{All eigenspaces in this proof are with respect to $h$.}
Lifting the eigenspace, $E_1 \neq 0$. By solubility, $E_0 \neq 0$. By Proposition~\ref{p:4:E0E1}, $\dim E_0 = \dim E_1 = 1$ and $E_2 = 0$, so $E_1$ is a $1$-dimensional, abelian subring. As in Proposition~\ref{p:4:b3}, which will be used heavily, we let $I = C_\fb^\circ(\fg/\fb)$.

\begin{step}\label{p:4:contradiction:st:hinE0}
We may assume $h \in E_0$.
\end{step}
\begin{proofclaim}
Consider the action of $E_0$ on $\fg/\fb$.

Suppose $E_0$ centralises $\fg/\fb$, viz.~$E_0 \leq I$. By Proposition~\ref{p:4:E0E1}, one has $\dim E_0 = \dim E_1 = 1$. If $E_1 \leq \fb$, then $I = E_0 + E_1$ where neither is trivial, against Proposition~\ref{p:4:b3}~\ref{p:4:b3:k-l}. Hence $E_1 \not \leq \fb$ and $(E_1\cap \fb)^\circ = 0$. However $E_0 \leq I$, so $[E_0, E_1] \leq (\fb\cap E_1)^\circ = 0$ and $E_0$ centralises $E_1$. Then $E_1 \leq C_\fg^\circ(E_0)$. The latter contains $I$; by Proposition~\ref{p:4:b3}~\ref{p:4:b3:st:buniqueonI}, it equals either $I$ or $\fb$. This is a contradiction in either case.

Therefore $E_0 \not\leq I$, and $\fb = E_0 + I$. Write $h = e_0 + i$ in obvious notation. Then $e_0$ also acts on $\fg/\fb$ like the identity; hence $E_1(e_0) \neq 0$. Moreover $E_0 = E_0(h) \leq E_0(e_0)$ by abelianity. By Proposition~\ref{p:4:E0E1}, one has $\dim E_0(e_0) = 1$, so $E_0(e_0) = E_0(h)$ contains $e_0$. Up to considering $e_0$ we are done.
\end{proofclaim}

\begin{step}\label{p:4:contradiction:st:anotleqb}
$E_1 \not\leq \fb$.
\end{step}
\begin{proofclaim}
Suppose $E_1 \leq \fb$. Let $f = \ad_h - 1$, an additive endomorphism of $\fg$. Let $K_1 = \ker^\circ f = E_1$ and $K_2 = \ker^\circ(f^2) \geq K_1$. By definition of $h$, one has $\im f \leq \fb$. Since $\dim \ker f = 1$, one has $\im f = \fb$, viz.~$f\colon \fg \twoheadrightarrow \fb$ is surjective. Since $K_1 = E_1 \leq \fb$, the restriction of $f\colon K_2 \to K_1$ remains surjective. Therefore $\dim K_2 = 2 \dim K_1 = 2$.

We claim $K_2 = I$. Let $k_2 \in K_2$, $k_1 = f(k_2) \in K_1$, and $e \in E_1$. Then:
\[[h, [e, k_2]] = [he, k_2] + [e, hk_2] = [e, k_2] + [e, k_2 + k_1] = 2[e, k_2] + [e, k_1].\]
However $e, k_1 \in E_1 = K_1$ while $E_2 = 0$, so $[e, k_1] = 0$. There remains $[e, k_2] \in E_2 = 0$, whence $K_2 \leq C_\fg^\circ(E_1)$.
Now $E_1 \sqsubset I$ so $I \leq C_\fg^\circ(E_1)$. If the latter has dimension $3$, then it is a $3$-dimensional Borel subring normalised by $h$. By Lemma~\ref{l:selfnormalisation}, $h \in C_\fg^\circ(E_1)$, a contradiction. Therefore $C_\fg^\circ(E_1) = I = K_2$. In particular, $B_{k_2} \leq \fb$.

In the notation above, also let $x \in \fg$. By definition of $h$, there is $b \in \fb$ with $[h, x] = x + b$. Therefore:
\[[h, [k_2, x]] = [hk_2, x] + [k_2, hx] = [k_2 + k_1, x] + [k_2, x + b] = 2[k_2, x] + [k_1, x] + [k_2, b].\]
The left-hand is in $\fb'$. So is $[k_2, b]$.
If $E_{-1} = 0$ then $[h, B_{k_1}] = B_{k_1}$ so $B_{k_1} \leq \fb'$. In that case, $[k_2, x] \in \fb' \leq I$, and therefore $K_2 = I$ is an ideal of $\fg$, a contradiction. So $E_{-1} \neq 0$. Let $\fh = E_{-1} + E_0 + E_1$, which must be soluble. Hence $E_{-1} \leq C_\fg^\circ(E_1) = I$ and $\fb = E_0 + E_{-1} + E_1$. But then for $e_{-1} \in E_{-1} \leq I = K_2$, one has $f(e_{-1}) = - 2 e_{-1} \in E_1$, whence $e_{-1} = 0$, a contradiction.

This proves $E_1 \not\leq \fb$.
\end{proofclaim}

Let $\fa = E_1$; by Step~\ref{p:4:contradiction:st:anotleqb} one has $\fg = \fa + \fb$. Let $\fb_2 = \fa + E_0$, a soluble subring. By Step~\ref{p:4:contradiction:st:hinE0} we could assume $h \in E_0$; so one has $\fb' = [E_0, \fa] = \fa$.
Let $\fb_2 < V \leq \fg$ be a $\fb_2$-irreducible module; $V$ need not be a subring. Linearising in dimension $1$ or by Corollary~\ref{c:Lie2}, $[\fb'_2, V] \leq \fb_2$.

Let $Y = (I \cap V)^\circ$, an $E_0$-module. Since $I$ is abelian by Proposition~\ref{p:4:b3}~\ref{p:4:b3:st:Iabelian}, $Y$ is actually a subring. A priori, $0 < Y \leq I$.
Notice $[\fa, Y] \leq (\fb\cap \fb_2)^\circ = E_0$.
We also claim $[E_0, Y] = Y$ since otherwise $h$ centralises $Y$, so $Y \leq E_0$ and equality holds, forcing $h \in I$, a contradiction.

If $\dim Y = 2$ then $Y = I$. Then for $a \in \fa$, the map $\ad_a$ takes $I$ to $E_0$. But by Proposition~\ref{p:4:b3}~\ref{p:4:b3:st:Z} one has $\dim Z(\fb) = 0$, so there is $i \in I\setminus Z(\fb)$ with $[a, i] = 0$. By~Proposition~\ref{p:4:b3}~\ref{p:4:b3:st:Iabelian}, $a$ normalises $C_i = I$ and therefore $a$ normalises $N_\fg^\circ(I) = \fb$. Hence $a \in \fb$ by Lemma~\ref{l:selfnormalisation}. This proves $\fa \leq \fb$, against Step~\ref{p:4:contradiction:st:anotleqb}.
Therefore $\dim Y = 1$.

Let $\fh = \fa + E_0 + Y$, a $3$-dimensional subring. By Proposition~\ref{p:4:minimalsimple}, $\fh$ is soluble. Since $[E_0, Y] = Y$ and $[E_0, \fa] = \fa$, solubility implies $[\fa, Y] = 0$. So $\fa \leq C_\fg^\circ(Y) = I \leq \fb$ by Proposition~\ref{p:4:b3}~\ref{p:4:b3:st:buniqueonI}, a final contradiction.
\end{proof}

It would be interesting to push this line further and try to classify minimal simple (or better: $N_\circ^\circ$-)Lie rings of finite Morley rank. Care will be needed since there is no clear indication what unipotence theory will become in this context, and since the Witt algebra is $N_\circ^\circ$.

\section{Appendix: questions}\label{S:questions}

We list a number of questions, some already in \cite{DRegard}. 

\subsection*{Extending the result}

\begin{question}
What happens to our theorem in characteristic $3$?
\end{question}

This is likely to be a question for pure algebraists, and a challenging one. Fact~\ref{f:1} is open in characteristic $3$.

\subsection*{Abstract Lie rings and their representations}


\begin{question}
State and prove an analogue of the `Borel-Tits' theorem \cite{BTHomomorphismes}, viz.~describe automorphisms of Lie rings of Lie-Chevalley type. 
\end{question}

One could start or content oneself over algebraically closed fields (see the elegant model-theoretic proof for groups in \cite{PMessieurs}), and even ask the same in Lie-Cartan type.

\begin{question}\label{q:CPT}
Devise identification methods à la Curtis-Phan-Tits \cite{TCurtis, GDevelopments} for Lie rings of Lie-Chevalley type.
\end{question}

\begin{question}
Study Lie ring representations of Lie rings of Lie-Chevalley type $\fG_\Phi(\bK)$. (Started in \cite{DSymmetric1, DLocally}. Also see Question~\ref{q:Sasha}.)
\end{question}

One could also ask about Lie rings of Lie-Chevalley type (or even Lie-Cartan type) up to elementary equivalence, à la Malcev \cite{MElementary, BElementary}. Such questions usually rely on retrieving the base field, which should be easier here than in groups.
For more on elementary properties of Chevalley \emph{groups}, see \cite{BMElementary}.

\subsection*{Lie modules of finite Morley rank}

In the proof of Theorem~\ref{t:4} (more specifically, when analysing $3$-dimensional Borel subrings), results like relevant analogues of Maschke's Theorem \cite[Proposition~2.13]{TStructure} or Lie's Theorem (in characteristic larger than Morley rank, generalising Corollary~\ref{c:Lie2}) would come handy.

\begin{question}
Develop basic tools for Lie modules of finite Morley rank.
\end{question}

\begin{question}\label{q:Sasha}
Let $\bK$ be an algebraically closed field of characteristic $p > 0$ and $\fg$ be a finite-dimensional, simple $\bK$-lie algebra. Assume that $\fg$ acts definably and irreducibly on an elementary abelian $p$-group $V$ of finite Morley rank.
Prove that:
\begin{enumerate}[label=(\roman*)]
\item\label{q:Sasha:i:linear}
$V$ has a structure of a finite dimensional $\bK$-vector space compatible with the action of $\fg$ [see \cite[Theorem~3]{BFinite}];
\item\label{q:Sasha:i:algebra}
$\fg$ is a $\bK$-Lie subalgebra of $\fgl_\bK(V)$. [More dubious.]
\end{enumerate}
\end{question}

The original result was for \emph{connected algebraic groups}, which correspond more to Lie rings of Lie-Chevalley type. And for that reason, we are more confident in~\ref{q:Sasha:i:linear} than in~\ref{q:Sasha:i:algebra}.

\begin{question}
Work in a theory of finite Morley rank.
Prove a `Steinberg tensor product theorem' \cite{SRepresentations} for definable $\fG_\Phi(\bK)$-modules. [See \cite[Theorem~3]{BFinite}.]
\end{question}

A more abstract direction makes no assumptions on $\fg$.

\begin{question}
Classify faithful, irreducible $\fg$-modules of Morley rank $2$. [See~\cite{DActions}.]
\end{question}

\subsection*{Abstract Lie rings of finite Morley rank}

There is of course the question of field interpretation in a non-abelian nilpotent Lie ring, where unreasonable hopes are shattered by \cite{BNew}.

Returning to the Reineke phenomenon, some topics are quite unclear.

\begin{question}\leavevmode
\begin{itemize}
\item
Is there a more model-theoretic proof of Fact~\ref{f:1} (viz.~one \emph{not} using Block-Premet-Strade-Wilson)?
\item
Let $\fg > 0$ be a connected Lie ring of finite Morley rank and $x \in \fg$. Is $C_x$ infinite? [See~\cite[Proposition~1.1]{BBCInvolutions}.]
\item
When is $x \in C_x$? 
[One should be careful with such questions, the Baudisch algebras could be counter-examples.]
\end{itemize}
\end{question}



At the soluble level, we have reasonable expectations.

\begin{question}
Develop a theory of soluble Lie rings of finite Morley rank: Fitting theory and Cartan subring theory. [See \cite{FAnormaux}.]
\end{question}

Of course we do not aim for `conjugacy' (if at all meaningful since there is no group around; fails anyway in Witt's algebra)---but existence and self-normalisation. Notice that existence could be hoped without the solubility assumption [see~\cite{FJExistence}], though we currently see no way to attack this.

Then of course, what matters is the simple case.
%

Let $\bK$ be an algebraically closed field of characteristic $p$. Then Witt's algebra $W_\bK(1; \underline{1})$ is a $p$-dimensional, simple Lie algebra over $\bK$, hence definable in the pure field.
Moreover, it has a subalgebra of codimension $1$.
Mind the assumption on the characteristic in the following.

\begin{question}
Let $\fg$ be a simple Lie ring of finite Morley rank $d$ and characteristic $p > d$. Suppose $\fg$ has a definable subring of corank $1$. Then $\fg \simeq \fsl_2(\bK)$. [See~\cite{HAlmost}.]
\end{question}

Last but not least in this direction, is of course the question of finding a strategy for the `$\log \CZ$' conjecture. For $\CZ$ itself, the positive solution in so-called `even type' \cite{ABCSimple} was a form of the classification of the finite simple groups (\textsc{cfsg}) seen from the \emph{reducing} lens of model-theoretic arguments for the tame infinite. At the high end of the classification, viz.~in high Lie rank, are `generic identification results'.

\begin{question}
Devise generic identification methods for simple Lie rings of finite Morley rank. [See Question~\ref{q:CPT}, see~\cite{BBGeneric}.]
\end{question}

As said in the introduction, we are curious whether model theory may provide a reduced sketch of the Block-Premet-Strade-Wilson theorem, like it did for the \textsc{cfsg}.

\subsection*{Model theory}

We see two questions here; one aims at broadening the model-theoretic frame of study of Lie rings, and the other at tightening connections between groups and Lie rings, presumably under new assumptions.

First, it is unclear to us whether the full strength of finite Morley rank was used; dimensionality could suffice, but this may depend on the behaviour of differential fields, if any. (It is unclear to us whether non-trivial differential fields would give counter-examples to the dimensional version of our results. Reasoning by loose analogy, fields with non-minimal $\mathbb{G}_m$ have not provided counter-examples to the Cherlin-Zilber conjecture, at least not yet.)

\begin{question}
Determine in which model-theoretic setting the present paper took place.
\end{question}

Second, by Lie-Chevalley correspondence, we mean the ability to attach to some groups a Lie ring. We do not believe in one for abstract groups of finite Morley rank, as categorical nature allows no infinitesimal methods. (This contrasts sharply with $o$-minimal nature, where infinitesimals are provided by elementary extensions; in the tame \emph{ordered} case, model-theoretic functoriality proves geometrically sufficient.)
One could add a stronger model-theoretic assumption with geometric flavour.

\begin{question}\label{q:LieChevalleyZariski}
Is there a Lie-Chevalley correspondence for groups definable \emph{in Zariski geometries} \cite{HZZariski}?
\end{question}

To our understanding, the Cherlin-Zilber conjecture is still open for groups definable in Zariski geometries. A positive answer to Question~\ref{q:LieChevalleyZariski} would certainly give a natural proof.

\subsection*{Acknowledgements}

We thank Gregory Cherlin for useful comments.

\paragraph{Very special thanks: the CIRM.}
The collaboration started virtually during a lockdown; but for serious research one must meet in person. This was made possible thanks to the \textsc{cirm-aims} residence programme (no longer restricted to host \textsc{aims} countries), allowing the second author to visit France for one week and the first author to visit Gabon for two, all in the Summer of 2023.
Both stays proved highly productive and pleasant.
Our very special thanks to Olivia Barbarroux, Pascal Hubert and Carolin Pfaff.

\printbibliography

\end{document}